\documentclass[11pt]{article}
\usepackage{amsfonts}
\usepackage{soul}
\usepackage{amssymb}
\usepackage{amsmath}
\usepackage{graphicx, subfigure}
\usepackage{color}
\usepackage{parskip}
\usepackage{graphicx, subfigure}
\newlength{\dinwidth}
\newlength{\dinmargin}
\newtheorem{definition}{Definition}
\newtheorem{theorem}{Theorem}

\newtheorem{proposition}{Proposition}

\newtheorem{lemma}{Lemma}
\newtheorem{example}{Example}

\begin{document}
\title{Polygons of Petrovi\'c and Fine, algebraic ODEs, and contemporary mathematics}

\author{Vladimir Dragovi\'c$^1$, Irina Goryuchkina$^2$}

\date{}

\maketitle

\footnotetext[1]{Department of Mathematical Sciences, The University
of Texas at Dallas, 800 West Campbell Road, Richardson TX 75080,
USA. Mathematical Institute SANU, Kneza Mihaila 36, 11000
Belgrade, Serbia.  E-mail: {\tt
Vladimir.Dragovic@utdallas.edu} 
}

\footnotetext[2]{Keldysh Institute of Applied Mathematics of RAS,  Moscow, Russia. E-mail: {\tt igoryuchkina@gmail.com }
}

\begin{abstract}

Here, we study the genesis and evolution of geometric ideas and techniques in investigations of movable singularities of algebraic ordinary differential equations. This leads us to the work of Mihailo Petrovi\'c on algebraic differential equations and in particular his geometric ideas captured in his polygon method from the last years of the XIXth century, which have been left completely unnoticed by the experts. This concept, also developed in a bit a different direction and independently by Henry Fine, generalizes the famous Newton-Puiseux polygonal method and applies to algebraic ODEs rather than algebraic equations. Although remarkable, the Petrovi\'c legacy has been practically neglected in the modern literature, while the situation is less severe in the case of results of Fine. Thus, we study the development of the ideas of Petrovi\'c and Fine and their places in contemporary mathematics.


\end{abstract}

\tableofcontents

\section{Introduction}

Here, we study the genesis and evolution of geometric ideas and techniques related to movable singularities of ordinary differential equations (ODEs). This leads us to the work of Mihailo Petrovi\'c on algebraic differential equations and in particular his geometric ideas captured in his polygonal method from the last years of the XIXth century, which have been left completely unnoticed by experts. This concept, also developed  independently by Henry Fine in a bit different direction, generalizes the famous Newton-Puiseux polygonal method (\cite{Newton},\cite{Cramer}, \cite{Puiseux}, see also contemporary sources like \cite{BK1986}, \cite{Ghys}. \cite{KK2012} and references therein)  and applies to algebraic ODEs rather than algebraic equations. Mihailo Petrovi\'c (1868-1943) was an extraordinary person and the leading Serbian mathematician of his time. His results are, despite their significance, practically unknown to  mathematicians nowadays. The situation is less severe with Fine's results. Thus, we emphasize here the development of the ideas of Petrovi\'c and Fine from the point of view of contemporary mathematics.

In its essence, this is a story about two outstanding individuals, Henry Fine and Mihailo Petrovi\'c, and their fundamental contributions to the development
from the point of view of modern mathematics of an important branch of analytic theory of differential equations. Along with doing science,  both made transformational efforts in elevating  of the mathematical research in their native countries to a remarkable new level. At the same time both left the deepest trace in the development of their own academic institutions at the moments of their transformation from a local college to a renewed university: Fine to US mathematics and Princeton University and Petrovi\'c
to Serbian mathematics and the University of Belgrade. The list of striking similarities between the two scientists is not even closely exhausted here. Both mathematicians were sons of a theologian. Both were in love with their own two beautiful rivers. Both actively enjoyed music playing their favorite instruments. Both went abroad, to  Western Europe, to get top mathematicians of their time as mentors to do their PhD theses. And both had a state official of the highest rank in their native countries as the closest friend. These friendships heavily shaped their lives. Fine published five scientific papers and had no known students. Petrovi\'c published more than 300 papers and had more than 800 scientific decedents.
It seems surprisingly that the two did not know each other and did not know about each other's work.

\subsection{Mihailo Petrovi\'c}

Mihailo Petrovi\'c defended his PhD thesis at the Sorbonne in 1894. His advisors were Charles Hermite and Emil Picard. During his studies in Paris, he learned a lot from Paul Painlev\'e and Julius Tannery and they become friends later. He was the founding father of the modern Serbian mathematical school and one of the first eight full professors of the newly formed University of Belgrade. At the same time he was a world traveler who visited both North and South Poles and a very talented travel writer.  Petrovi\'c also received  many awards as an inventor, including the gold medal at the World Exposition in Paris 1900. His nickname Alas referred to his second profession, a fisherman on the Sava and the Danube rivers. He was even more proud of his fisherman's achievements than any others, and Alas became an integral part of his full name. He played violin and founded a music band called Suz. Professor Petrovi\'c was a member of the Serbian Academy of Sciences and Arts, a corresponding member of the Yugoslav Academy of Sciences and Arts, and a foreign member of several other academies. He was invited speaker at five international congresses of mathematicians,  Rome 1908, Cambridge 1912, Toronto 1924, Bologna 1928, and Zurich 1932. Professor Petrovi\'c was also the personal teacher  and mentor of the then crown prince George of Serbia, with whom he remained friends for life. Their friendship became even stronger after George abdicated in favor of his younger brother and later king, Alexander I of Yugoslavia (see \cite{Karadjordjevic}). Professor Petrovi\'c's legacy included 11 PhD students and almost 900 PhD students of his PhD students and their students. He was a veteran of Balkan wars and  the First World War. The Serbian and later Yugoslav army used his cryptography works for many years.  In 1941, when the Second World War arrived in his country, he got mobilized as a reserve officer. After the Axis powers occupied his country, Petrovi\'c ended as a prisoner of war in Nuremberg, at the age of 73. The former crown prince George made a plea to his aunt, the queen Elena of Italy, based on Petrovi\'c's illness. Thus, Petrovi\'c got released from the prisoners camp and soon died in Belgrade, his place of birth.

A street in the downtown, an elementary school, a high-school, and a fish-restaurant in Belgrade are named after Mihailo Petrovi\'c Alas.

The main goal of this paper is to introduce modern readers to the results of Mihailo Petrovi\'c and to relate them to the results of modern theory of analytic differential equations. The main source of Petrovi\'c's results for us is his doctoral dissertation \cite{Petrovich1}. It was written in French in 1894. (It was reprinted along with a translation in Serbian, edited by Academician Bogoljub Stankovi\'c in Volume 1 of \cite{PetrovicCW} in 1999.)  The dissertation consists of two parts. The first part of the thesis mostly presents results related to the first order algebraic ODEs, while the second part is related to algebraic ODEs of higher orders.  At the beginning of the first part of the thesis, Petrovi\'c introduces a new method. We are going to call it {\it the method of Petrovi\'c polygons}. This method is suitable to study analytic properties of solutions of algebraic ODEs in a neighborhood of the  nonsingular  points of the equation. The method of Petrovi\'c polygons is a certain modification of the method of Newton-Puiseux, applicable to the study of solutions to algebraic equations. Further on, Petrovi\'c applies his method to study zeros and singularities of algebraic ODEs of the first order. He formulates and proves necessary and sufficient conditions for the { non-existance of movable zeros and poles}  of solutions (see Theorem \ref{th:1} below). Contrary to the  necessary and sufficient conditions for the { non-existance} of movable critical points of solutions of algebraic ODEs of the celebrated L. Fuchs Theorem \cite{Fuchs} (provided below as Theorem \ref{th:Fuchs}), the conditions of Petrovi\'c's Theorem \ref{th:1} do not request either the computation of solutions of the discriminant equation or to have the equation resolved with respect to the derivative. The conditions of Petrovi\'c's Theorem \ref{th:1}  can be checked  easily and effectively by a simple constructing of a geometric figure corresponding to the given equation.  The first part of the dissertation also contains the theorems which provide a classification of rational, first order ODEs explicitly resolved with respect to the derivative which have uniform (single-valued) solutions (see Petrovi\'c's Theorem  \ref{th:4} below). Later on, these results of Petrovi\'c were essentially improved by J. Malmquist \cite{Malquist} (see Theorem \ref{th:Malmquist} below). In addition, in the first part of the thesis, the class of binomial ODEs of the first order is studied and the equations with solutions without movable singular points are described. Also, Petrovi\'c characterizes those binomial equations which possess uniform (single-valued) solutions. The results of Petrovi\'c are very similar to those obtained by K. Yosida more than 30 years later  \cite{Yoshida} (see Theorem \ref{th:Yoshida} below). The second part of the dissertation is devoted to the applications of the polygon method in the study of zeros and singularities of the algebraic higher order ODEs. We will present some of the results related to higher order ODEs a bit later  (see Section \ref{sec:higherorder} and Theorems \ref{th:Acta1} and \ref{th:Acta2}).

Mihailo Petrovi\'c is one of the most respected and influential mathematicians in Serbia. Petrovi\'c's collected works in 15 volumes were published in 1999 \cite{PetrovicCW}. The year 2018 was the Year of Mihailo Petrovi\'c in Serbia on the occasion of his 150th anniversary. A monograph was published by the Serbian Academy of Sciences and Arts to celebrate his life and scientific results (see \cite{Monografija}). Nevertheless, it happened that none of his students and followers in
Serbia continued to develop further the  geometric ideas from his PhD thesis. The first two out of eleven of his direct PhD students, Dr Mladen Beri\'c and Dr Sima  Markovi\'c
were involved in the topics closely related to their mentor's dissertation, see
 \cite{Beric} and \cite{Marcovic}.  However, the life produced unexpected turns: Dr Beri\'c was forced to leave the University of Belgrade at the beginning of 1920s due to the issues related to his personal life. On the other hand, Dr Markovi\'c became the first secretary of the Yugoslav Communist Party. When the Communist Party was prohibited by law in 1920, Markovi\'c lost his position at the University of Belgrade immediately after. (Later on he came into a dispute with the Party. Finally he was executed after a quick trial in Moscow in 1939, to be rehabilitated in 1958.)   These extraordinary circumstances can at least partially explain the lack of continuity within the school founded by Petrovi\'c in the field of analytic theory of differential equations, although continuity exited in many other directions pursued by later students of Petrovi\'c, like Pejovi\'c, Mitrinovi\'c, Ka\v sanin, or Karamata. Serbian mathematicians who have been active nowadays in the field of analytic theory of differential equations (see, for example, \cite{DS2019},\cite{DGS2018}, \cite{JR1}, \cite{JR2}) neither methodologically nor according to the mathematical genealogy belong to the Petrovi\'c school.
 Certainly, some of the Petrovi\'c's results in that field were quite well known at the beginning of the XXth century. Nevertheless, neither Golubev, who extensively used some other results from Petrovi\'c's thesis in his famous book
 \cite{Golubev}, nor any other mathematician who used later analogous geometric methods in the study of the solutions of the algebraic ODEs, ever quoted Petrovi\'c's foundational results in this field (see  \cite{Cano2}, \cite{GrigorievZinger}, \cite{Bruno}).

\subsection{Henry Fine}

It should be pointed out that, a couple of years prior to Petrovi\'c, the American mathematician Henry Fine invented another modification of the Newton-Puiseux method for studying the formal solutions of algebraic ODEs \cite{Fine}. Let us notice that although the Fine construction was similar to the one of Petrovi\'c, they were not identical and the questions they were considering were very much different. As we have said, it seems that Petrovi\'c and Fine didn't know about each others results. At the end of the XXth century, the Fine method was developed further by  J. Cano \cite{Cano2}, \cite{Cano1}. As of today, contemporary methods based on different modification of the Newton-Puiseux polygonal method allow wide classes of formal solutions to be computed for analytic differential equations and their systems \cite{Bruno}, and to prove their convergence and analysis of the rate of growth of terms of formal series  \cite{Cano2}, \cite{Malgrange}, \cite{Ramis}, \cite{Sibuya}, which is needed to know in the selection of  summation methods.

Henry Burchard Fine (1858-1928) was a dean of faculty and the first and only dean of the departments of science at Princeton. He was one of a few who did most to help Princeton develop from a college into a university. He made Princeton a leading center for mathematics and fostered the growth of creative work in other branches of science as well. Professor Oswald Veblen, in his memorial article \cite{Veblen} described Fine's contribution on the nationwide scale in his opening sentence by saying that ``Dean Fine was one of the group of men who  carried American mathematics forward from a state of approximate nullity to one verging on parity with the European nations."

Fine grew up between two major rivers, the St. Lawrence and the Mississippi river and was always astonished by them. ``He played the flute in the college orchestra, rowed on one of the crews, and served for three years as an editor of the Princetonian, where he began a life-long friendship with Woodrow Wilson (in) 1879, whom he succeeded as managing editor." \cite{Leitch}

He traveled to Leipzig, Germany, in 1884 to  attended lectures by Klein, Carl Neumann and others and to prepare his doctoral dissertation. The topic was suggested by Study, and approved by Klein. Fine defended his dissertation ``On the singularities of curves of double curvature" in May 1885 at  the University of Leipzig. Upon his return from Germany,
Fine was appointed assistant professor at Princeton. Despite great promise as a research mathematician, Fine moved very soon into other areas of academic life. He mainly devoted his time to teaching, administration, and the logical exposition of elementary mathematics.

His first research paper came out of his thesis, had the same title ``On the singularities of curves of double curvature" and appeared in the American Journal of Mathematics in 1886. In the following year he published a generalization of these results to $n$ dimensions in the same journal. Two further papers ``On the functions defined by differential equations with an extension of the Puiseux polygon construction to these equations", and ``Singular solutions of ordinary differential equations" from 1889 and 1890 respectively, are of the utmost importance for our current presentation (see \cite{Fine}, \cite{Fine1}). His last research publication  appeared in 1917: ``On Newton's method of approximation".

Fine was one of the founders of the American Mathematical Society. He served as the AMS president in 1911 and 1912.

Fine wrote several books on elementary mathematics, including ``Number system of algebra treated theoretically and historically", ``A college algebra", ``Coordinate geometry", and ``Calculus".

We quote a few more very illustrative parts from \cite{Leitch}:
``In 1903, shortly after he became president of the University, Wilson appointed Fine dean of the faculty, and Fine's energies were thereafter devoted chiefly to university administration. He worked shoulder to shoulder with Wilson in improving the curriculum and strengthening the faculty, and bore the onus of the student dismissals made inevitable by the raising of academic standards. In the controversies over the quad plan and the graduate college, Fine supported Wilson completely. After Wilson resigned to run for governor of New Jersey in 1910, Fine, as dean of the faculty, carried the chief burden of the university administration during an interregnum of two years; and when the trustees elected John Grier Hibben as Wilson's successor, Fine, who many had thought would receive the election, magnanimously pledged Hibben his wholehearted support. `He was singularly free from petty prejudices and always had the courage of his convictions,' Hibben later recalled. `Every word and act was absolutely in character, and he was completely dependable in every emergency.'...After his election as president of the United States, Wilson urged Fine to accept appointment as Ambassador to Germany and later as a member of the Federal Reserve Board, but Fine declined both appointments, saying quite simply that he preferred to remain at Princeton as a professor of mathematics. Fine also declined a call to the presidency of Johns Hopkins University and several to the presidency of Massachusetts Institute of Technology....In the summer of 1928, he went to Europe, where he revisited old scenes and old friends, and recovered to some extent, in the distractions of travel, from the sorrows he had suffered in the recent death of his wife and the earlier deaths of two of his three children. Professor Veblen who talked with him soon after his return, reported later that Fine `spoke with humorous appreciation of the change he had observed in the attitude of European mathematicians toward their American colleagues and with pride of the esteem in which he had found his own department to be held.'

Tall and erect, Dean Fine was a familiar figure on his bicycle, which he rode to and from classes and used for long rides in the country. While riding his bicycle on the way to visit his brother at the Princeton Preparatory School late one December afternoon, he was struck from behind by an automobile whose driver had failed to see him in the uncertain light of dusk. He died the next morning..."

The Mathematics Department of the Princeton University is housed in the Fine Hall, the building named after its first chairman.

\section{Historic context: the 1880's}
Across the entire XIXth century there was a significant and constant attention to the study of analytic differential equations and their classification. The French school with Cauchy, Liouville, Picard, Hermite, Briot, Bouquet,  Poincar\'e, and Painlev\'e gave tremendous contributions. Of course, among those who gave a key contribution to the analytic theory of differential equations were scientists from other countries as well, such as Gauss, Riemann, Weierstrass, L. Fuchs, and Kowalevski. Nevertheless, the center of the attention to analytic theory of differential equations was indeed in France. And that was definitely the case in 1880's, when French mathematics contributed greatly to further development of complex analysis and to applications
of its methods to the study of differential equations. So called algebraic differential equations and systems of such equations occupied
a special place. Let us recall that an algebraic ODE has the following form
\begin{equation}
\label{eq1}
f(x,y,y',\dots,y^{(n)})=0,
\end{equation}
\begin{equation}
\label{eq2}
f=\sum\limits_{i=1}^{j}\varphi_i(x)\,y^{m_{0i}}y'^{m_{1i}}\cdots y^{(n) m_{ni}},
\end{equation}
$j\in\mathbb{N},$ where $m_{0i},\dots,m_{ni}\in\mathbb{Z}_+,$ and $\varphi_i$ -- are algebraic or analytic functions.

Let us also recall the main notions of the analytic theory of differential equations (see for example \cite{Golubev}), which are necessary to formulate the results of Petrovi\'c and his predecessors.

  The points where a given solution of the equation \eqref{eq1} is not analytic or it is not defined are called {\it the singular points of the solution}. Simple examples of singular points of a solution are its poles, i.e. the points such that in their punctured neighborhoods the solution is presented by Laurent series with finite principal parts. Paul Painlev\'e suggested to classifying isolated singular points of a function according to the number of values it takes while going around the singular points. This led to the division of singular points on {\it critical} and {\it noncritical} points. If the function changes its values while going around a singular point, the singular point is {\it critical}. Examples: the point $x=0$ is a critical point for the function  $\sqrt{x}$ and for the function $\ln{x}$. If, on the contrary, a function does not change its value while going around a singular point, the singular point is called {\it noncritical}. Examples: the point $x=0$ is a noncritical singular point for functions $1/x$ and $e^{1/x}$.
If going around a critical singular point $x=a$, the function takes a finitely many values and has a limit in that point (as $x \rightarrow a$ inside any sector with the vertex in $a$ and with a finite angle), then such a  point is either an algebraic critical point or a critical pole. It is called {\it an algebraic critical point} if in its neighborhood the function has an expansion of the form
 $$y=c_0+c_1(x-a)^{1/s}+c_2(x-a)^{2/s}+\dots$$
 (i.e. the limit is finite).
 The point $x=a$ is called \textit{a critical pole}, if in a punctured neighborhood of the point, the function has an expansion of the form: $$y=c_{-m}(x-a)^{-m/s}+\dots+c_{-1}(x-a)^{-1/s}+c_0+c_1(x-a)^{1/s}+\dots,\qquad {s\geqslant 2} $$
 (i.e. the limit is infinite).
 The critical algebraic points, poles, and critical poles form the set of {\it algebraic singular points}. {

 In turn, as seen from the above definitions, {\it nonalgebraic critical points} are of two types: the first ones are such that the function takes an infinitely many values if going around a critical point ({\it transcendental points}; for example, the point $x=0$ for the function $\ln{x}$), the second ones are those at which the function does not have a limit ({\it essential singular points}; for example, the point $x=0$ for the function $\sin(1/\sqrt{x})$).

 Here we do not consider non-isolated singular points of functions. The question of the classification of such singular points is more delicate. Some examples of the application of the above classification of isolated singular points to studying singular points forming {\it singular curves} can be found in \cite{Golubev}, Ch. I, \S 7.}

It is well known that for the linear ODEs the singular points of the solutions form a subset of the set of singular points of the coefficients of the  equations. For nonlinear ODEs, the points where the coefficients  $\varphi_i(x)$ of the equation \eqref{eq1} are zero or undefined, as a rule, are singular points of its solutions. We will call such points {\it the singular points of the equation}. But it is important to stress that not only the singular points of a nonlinear ODE can be singular points of its solutions. This property of nonlinear ODEs led L. Fuchs to divide the singular points of the solutions of a nonlinear ODE into movable and fixed singular points. \textit{A fixed  singular point of a solution} of the equation \eqref{eq1} is a such its singular point whose position does no depend on initial data, which determine the solution, i.e. such a singular point is a common singular point for  $n$-parametric family of solutions, or in yet another words, the common singular point for the general solution (the general solution is also called the integral) of the equation.  Example: the point $x=0$ is a fixed singular point of the general solution $y=(C-\ln x)^2$, where $C$ is an arbitrary constant corresponding to the initial data, for the equation $x^2y'^2-4y=0$; the point $x=0$ is a singular point of the last equation.  {\it A movable singular point of a solution} is such an its singular point, whose position depends of the constants of integration. For example, $x=-C$\, where $C$ is an arbitrary constant, is  a movable singular point of the solution  $y=1/(x+C)$  of the equation $y'+y^2=0$. The last equation does not have singular points.

For nonlinear ODEs of the first order, {like for linear ODEs of any order}, it can be shown that fixed singular points are always singular points of the equation, see for example {Ch. I, \S 8} of \cite{Golubev}.

 As it has already been said, Petrovi\'c was concerned about the conditions on algebraic ODEs, under which their solutions would have or would not have movable singular points. He was also studying the nature of solutions, whether they are single-valued, rational, or elliptic functions.  The importance of movable critical points lies in the fact that their presence prevents the possibility for the given differential equation to construct a unique Riemann surface, which could serve as the domain for all the solutions of the given equation as uniform (single-valued) functions. Now we are going to list the results which had been widely known at the moment when Petrovi\'c arrived in Paris in 1889 to pursue his graduate education and which had systematically been used in his doctoral dissertation.

\begin{theorem}[Small Picard Theorem, \cite{PicardT1}, \cite{PicardT2}, \cite{Encyc}]
Outside the image of any nonconstant entire function there could be at most one complex number.

\end{theorem}

\begin{theorem}[Big Picard Theorem, meromorphic version, \cite{PicardT1}, \cite{PicardT2}, \cite{Encyc}]
 Let
$\mathcal{M}$ be a Riemann surface and $\mathcal{S}$ the Riemann sphere, and $w\in \mathcal {M}$ an arbitrary point.
Let $f:\mathcal{M}\setminus\{w\} \to \mathcal{S}$ be a holomorphic function with an essential singularity at the point $w$.
Then, all points on the sphere $\mathcal{S}$   except at most two have infinitely many inverse images.

\end{theorem}

 Let us mention some important results about the first order algebraic differential equations, used frequently by Petrovi\'c.

\begin{theorem}[Hermite, \cite{Golubev}, \cite{Hermite}]\label{th:Hermite} Given a { polynomial} $P$. If solutions of the equation  $P(y,y')=0$ do  not have movable critical points, then the genus of the curve $P(x, y)=0$ is equal either 0 or 1. In such a  case the solutions are either rational functions or could be rationally expressed through exponential and elliptic functions.
\end{theorem}

\begin{theorem}[H. Poincar\'e, L. Fuchs, 1884-5, \cite{Fuchs}, \cite{Golubev}, \cite{Poincare}]
Algebraic differential equations of the first order without movable critical points reduce to linear, Riccati (see \eqref{eq:re}), or Weierstrass equations (see \eqref{eq:we}).

\end{theorem}

In 1884 Lazarus Fuchs provided necessary conditions of absence of movable critical points of solutions of algebraic ODEs of the first order { \cite{Fuchs}. In fact, his theorem provided also sufficient conditions of absence of movable {\it algebraic} critical points of  solutions.  But a 1887 result of Painlev\'e allowed to strengthen the assertion of Fuchs theorem up to sufficient conditions of absence of {\it any} movable critical points. The result of Painlev\'e is also going to be stated a  paragraphs below as Theorem \ref{th:Painleve}. Now we are going to formulate the corresponding Theorem of Fuchs which is going to include the sufficient part as well.}

 Let
\begin{equation}\label{eq:thFuchs}
 F(x,y,y')=A_0(x,y) (y')^s+A_{1}(x,y)(y')^{s-1}+\dots+A_s(x,y),
 \end{equation}
 where { $A_0
 ,\dots, A_s$ are  polynomials in  $y$} and analytic in $x$.
{\it The discriminant equation} $D(x,y)=0$ is the result of the elimination of $y'$ from the equations $F(x,y,y')=0$ and $\displaystyle \frac{\partial F(x,y,y')}{\partial y'}=0$.

 \bigskip
 \begin{theorem}[L. Fuchs, 1884-5, \cite{Fuchs}, \cite{Golubev}] \label{th:Fuchs} The solutions of the equation  $F(x,y,y')=0$, where $F$ is defined in \eqref{eq:thFuchs}, don't have movable critical points if and only if:
 \begin{itemize}
 \item The coefficient {$A_0(x,y)$} does not depend on $y$.
 \item The degree of each of the polynomials  $A_k(x,y)$ with respect to  $y$ does not exceed $2k$.
 \item The solutions $\phi(x)$ of the discriminant equation $D(x,y)=0$ have to be integrals of the given equation.
 \item { For each fixed $x_0$, the Puiseaux expansion of $y'$ with respect to $y$ in a neighbourhood of the point   $y_0=\phi(x_0)$  has the form
 $$y'= \phi'(x_0)+b_k(x_0)(y-y_0)^{k/m}+b_{k+1}(x_0)(y-y_0)^{(k+1)/m}+{\dots}$$ with $k\geqslant m-1.$}
 \end{itemize}
\end{theorem}

{  Let us give a proof of the ``only if'' part of Fuchs theorem (one can see a proof in \cite{Fuchs} or in \cite{Golubev}).

\bigskip
 First we suppose that $A_0(x,y)$ contains $y$. We take an arbitrary value $x_0$ which differs from singular points of the equation and a value $y_0$ that satisfies $A_0(x_0,y_0)=0.$ As known from the course of analytic theory of differential equations (see, for example, \S~1, Ch.~II in \cite{Golubev}) near such points $(x_0,y_0)$ the equation $F(x,y,y')=0$ is resolved with respect to $y'$
 in the form
\begin{equation}\label{ext1}
y'=\frac{1}{c_l(y-y_0)^l+c_{l+1}(y-y_0)^{l+1}+\dots}+\alpha,\;\; \mbox{if}\;\; D(x_0,y_0)\neq 0,
\end{equation}
where $c_j\in\mathbb{C}\{x-x_0\},$ $\alpha\in\mathbb{C},$  or in the form
\begin{equation}\label{ext2}
y'=\frac{1}{c_l(y-y_0)^{l/m}+c_{l+1}(y-y_0)^{(l+1)/m}+\dots}+\alpha,\;\; \mbox{if}\;\; D(x_0,y_0)= 0,
\end{equation}
where $c_j\in\mathbb{C}\{x-x_0\},$ $\alpha\in\mathbb{C}$.
In the equations \eqref{ext1} and \eqref{ext2} the coefficient $c_l(x_0)\neq 0$ due to the point $x=x_0$ is not a singular point of the equation.

We can apply the lemma below to the equation \eqref{ext1} and,
after the change of variable $(y-y_0)^{1/m}=u$,  to the equation \eqref{ext2}.

\begin{lemma}[Lemma about a critical point, \cite{Golubev}]\label{lemma:lemma}
Let  $f(x_0, y_0)=\infty$ and $1/f(x, y)$ be holomorphic in a neighborhood of the point $(x_0, y_0)$. Then  $x_0$ is a movable critical algebraic point of the equation $y'=f(x,y)$. More precisely, the integral determined by the initial data $(x_0,y_0)$ has an expansion $$y=y_0
+a_1(x-x_0)^{1/k}+a_2(x-x_0)^{2/k}+\dots$$ near the point $x_0$, where $k=1+{\rm ord}_{y=y_0}(1/f(x_0,y))$.
\end{lemma}

According to this lemma the both equations have movable critical points of their solutions. Hence if a solution of the equation { $F(x,y,y')=0$} has only fixed critical points then the coefficient $A_0(x,y)$ does not contain $y$.
We prove the necessity of the first condition of the Fuchs theorem.

\bigskip
Further we write the equation $F(x,y,y')=0$ in the form
\begin{equation}\label{eq:thFuchs2}
 (y')^s+A_{1}(x,y)(y')^{s-1}+\dots+A_s(x,y)=0,
 \end{equation}
 where $A_1
 ,\dots, A_s$ are  polynomials in  $y$ and analytic in $x$.

 Let $p_j$ be the degree of the polynomial $A_j$ with respect to $y$. The coefficients
 $$A_j(x,y)=A_j(x,1/w)=\frac{B_j(x,w)}{w^{p_j}},$$ where $B_j(x,w)$ is a polynomial with respect to $w$.
  So by means of the transformation $y=1/w$ the equation \eqref{eq:thFuchs2} reduces  to the equation
  \begin{equation}\label{eq:thFuchs3}
  (w')^s-\frac{B_1(x,w)}{w^{p_1-2}}(w')^{s-1}+\frac{B_2(x,w)}{w^{p_2-4}}(w')^{s-2}-{\dots}+(-1)^{s}\frac{B_s(x,w)}{w^{p_s-2s}}=0.
  \end{equation}
As a solution of the equation \eqref{eq:thFuchs3} has no movable critical points then this equation does not contain $w$ in the denominators of the coefficients $\displaystyle \frac{B_j(x,w)}{w^{p_j-2j}}$.
Thus we prove the necessity of the second condition of the Fuchs theorem, that the equation {$F(x,y,y')=0$}  must be of the form
\eqref{eq:thFuchs2} with the polynomials $A_j$ of degree $p_j\leqslant 2j.$

\bigskip
The equation \eqref{eq:thFuchs2} with the polynomials $A_j$ of degree $p_j\leqslant 2j$ can be resolved with respect to the $y'$ in the form
 \begin{equation}\label{eq:thFuchs4}
 y'= s_0+b_k(y-v)^{k/m}+b_{k+1}(y-v)^{(k+1)/m}+{\dots},
\end{equation}
where { $F(x,v,s_0)=0$}, $v$ is an arbitrary function with respect to $x$, $b_j$ are holomorphic near $x_0$, $k,\,m\in\mathbb{N},$ $b_k\neq 0.$
 If $D(x,v)\neq 0$ then  $m=1$. If $D(x,v)= 0$ then $m>1$.

  Making the change of variable $y-v=w^m$ (whence $y'=mw^{m-1}w'+v'$) in the equation \eqref{eq:thFuchs4} we get the equation
\begin{equation}\label{ext5}
  w'=\frac{s_0-v'+b_k w^k+b_{k+1} w^{k+1}+\dots}{mw^{m-1}}.
\end{equation}

First we consider the case $s_0\not\equiv v'$.
If $s_0(x_0)-v'(x_0)\neq 0$  then according to the Lemma 1 the equation \eqref{ext5} with the initial data $w(x_0)=0$ has the
solution $$w=a_1(x-x_0)^{1/m}+a_2(x-x_0)^{2/m}+\dots,\;\;a_1\neq 0.$$

If $s_0(x_0)=v'(x_0)$ (say, $s_0-v'=(x-x_0)^n\varphi(x),$ $\varphi(x_0)\neq 0$) then changing the variable $w=(x-x_0)^nu$ in the equation \eqref{ext5} and applying the ideas above to the modified equation we get that in this case the equation \eqref{ext5} has solutions $$w=a_l(x-x_0)^{l/m}+{a_{l+1}}(x-x_0)^{(l+1)/m}+\dots,\;\;a_l\neq 0.$$

Taking into account the change of variable $(y-v)^{1/m}=w$ and the equation \eqref{eq:thFuchs4} we  get
$$y'=s_0+b_k a_1^k(x-x_0)^{k/m}+\dots.$$
Integrating the last equality we find that if $k$ is not divided by $m$  then the solution $y$ has an algebraic critical point $x_0$.

So for the absence of critical movable points it is necessary for $k/m$ to be integer. In such a way we can prove that all other terms of the series in the right part of the equation  \eqref{eq:thFuchs4} contain $y-v$ in integer powers. Therefore in the case when the equation  \eqref{eq:thFuchs4} has fractional power exponents of $y-v$ {(that is, when $D(x,v)=0$)} for $s_0-v'\not\equiv 0$, the equation \eqref{eq:thFuchs2}  has integrals with movable algebraic critical points.
{ This means that a solution $v$ of the discriminant equation $D(x,v)=0$ should satisfy $v'\equiv s_0$, that is, $F(x,v,v')=0$.}

Hence we prove the necessity of the third condition of Fuchs theorem, namely, if the equation \eqref{eq:thFuchs2}  has  only fixed critical points then all solutions of the discriminant equation are solutions of this equation too.

In the last part of the proof we have that $s_0=v'$.
And the equation \eqref{ext5} has the form
\begin{equation}\label{ext7}
  w'=\frac{b_k w^k+b_{k+1} w^{k+1}+\dots}{mw^{m-1}}.
\end{equation}
If $k\geqslant m-1$ the equation \eqref{ext7} and also the equation \eqref{eq:thFuchs4} have holomorphic integral with initial data $w(x_0)=0$.

In the case $k<m-1$ the equation \eqref{ext7} can be written in the form
\begin{equation}\label{ext8}
  w'=\frac{b_k +b_{k+1} w+\dots}{mw^{m-k-1}}.
\end{equation}
 Its integral that vanishes { at} the point $x=x_0$ has the form
 $$w=a_1(x-x_0)^{1/(m-k)}+a_2(x-x_0)^{2/(m-k)}+\dots.$$
From this for the right part of the equation \eqref{eq:thFuchs4} we get an expansion
$$y'=s_0+b_k a_1^k(x-x_0)^{k/(m-k)}+\dots,\;\;b_k(x_0),\;a_1(x_0)\neq 0.$$
Integrating the last equation we get that in this case a solution of the equation \eqref{eq:thFuchs4} has an algebraic critical movable point. Thus, we prove the last condition of the Fuchs theorem: in the expansion of $y'$ in  fractional powers of $y-\phi(x)$, where $\phi(x)$ is a solution of the equation $D(x,y)=0$, one has
 $$y'= \phi'(x)+b_k(y-\phi(x))^{k/m}+b_{k+1}(y-\phi(x))^{(k+1)/m}+{\dots}$$
 with $k\geqslant m-1.$
The proof of the ``only if'' part of Fuchs theorem is finished.

}

Let us briefly mention basic biographic data about Lazarus Fuchs. He was born in 1833 in Moschina, the Grand Duchy of Posen of Kingdom of Prussia, nowadays Poland.
He worked on his PhD in Berlin University with Kummer as his advisor, from 1854 till 1858, when he defended a thesis on the lines of curvature on surfaces. In 1882 he returned to Berlin where he got a position of a full professor of the Berlin University. He was elected a member of Berlin Academy in 1884. From 1892 till his death, Fuchs served as the Editor-in-Chief of "Journal f\" ur die reine und angewandte Mathematik" (the Crelle's journal). He died in Berlin in 1902.

\bigskip
As it has been mentioned in the Introduction, P. Painlev\'e in his doctoral dissertation \cite{Painleve} in 1887 formulated and proved a remarkable result, which inspired many to further investigations of solutions of algebraic ODEs. Painlev\'e defended his dissertation on 10th of June of 1887 in front of the committee chaired by Hermite, with Appell and Picard as the members. The thesis was devoted to Picard, the mentor of Painlev\'e. The dissertation was published as a journal paper a year later, see \cite{Painleve0}. The central result of the dissertation is the following:

\begin{theorem}[P. Painlev\'e, 1887, \cite{Golubev}, \cite{Painleve}, \cite{Painleve0}]\label{th:Painleve}
Differential equations of the first order, algebraic with respect to the unknown function and its derivative, don't possess movable  non-algebraic singularities.

\end{theorem}

The proof uses the following result of Cauchy about the differential equation
$$
w'={ f( z,w)},
$$
where $f$ is holomorphic within the disks $|w|\le \rho >0$ and $|z|\le r>0$. For any $\rho_1, r_1$ such that  $0<\rho_1<\rho$ and $0<r_1<r$, there exists $M$
such that { $|f(z,w)|<M$} for all $w, z$ with $|w|=\rho_1$ and $|z|=r_1$. Then there exists { a solution} of the above differential equation $w_1=w_1(z)$, such that $w_1(0)=0$
and $w_1$ is holomorphic within the disk
\begin{equation}\label{eq:Cauchyestimate}
|z|<r_1{(1-\exp(-\rho_1/(2Mr_1)))}.
\end{equation}

The second ingredient of the proof is Painlev\'e's Theorem 3, from p. 36 {of the Painlev\'e dissertation \cite{Painleve}}, which states:
{\it If for all points $z_0$ of a region $S$, the function $f(z_0, u)$ has at most a discrete set of essential points $a_1, a_2, \dots, a_m,\dots$ with the coordinates depending on $z_0$ analytically, the root $u(z)$ of the equation  $f(z, u)=0$ { is} single-valued (or $n$-valued) in a  region $S'$ with the boundary $\sigma$, a subset of the interior of $S$, with a continuation across $\sigma$, with the poles (or critical algebraic points) as the singularities in $S'$. }

Then, on p. 41 {of \cite{Painleve}}, Painlev\'e formulates the following statement: {\it Given a
{ differential equation}
$$
\frac {du}{dz}=f(z, u),
$$
where $f$ is a single valued function when $z$ varies in $S$ and $u$ in the complex plain.  If for all points $z_0$ of a region $S$, the function $f(z_0, u)$ has at most a discrete set of points where it is not defined, with the coordinates depending on $z_0$ analytically, all integrals $u(z)$  are single-valued (or $n$-valued) in a  region $S'$ with the boundary $\sigma$, a subset of the interior of $S$, with a continuation across $\sigma$, with the poles (or critical algebraic points) as the singularities in $S'$.}

As an important class of examples, Painlev\'e gives the equations of the form:
\begin{equation}\label{eq:Painleverational}
\frac {dw}{dz}=\frac{P(w,z)}{Q(w,z)},
\end{equation}
where $P, Q$ are polynomials with respect to $w$. We will provide more details of the proof of the Painlev\'e theorem in this case, following
\cite{Golubev}. First, let $M_1$ denote the set of points which are singularities in $z$ of the coefficients of $P, Q$ as  polynomials in $w$ and the common zeros of all the coefficients of $Q$.  Let $M_2$ denote the set of $z$ points for which $P=0$ and $Q=0$ has common solutions. Similarly, changing the variable $w_1=1/w$ we get the equation
 $$
\frac {dw_1}{dz}=\frac{P_1(w_1,z)}{Q_1(w_1,z)}.
$$
Let $M_3$ denote the common zeros of $P_1(0,z)=0$ and $Q_1(0, z)=0$. Changing the variable $z_1=1/z$ we get the equation
 $$
\frac {dw}{dz_1}=\frac{P_2(w,z_1)}{Q_2(w,z_1)}.
$$
If $z_1=0$ contributes to one of the sets $M_1$, $M_2$, or $M_3$ we add it as well and then the union of these sets we denote as $\mathcal M$.
Let us assume that a point $z_0$ outside $\mathcal M$ is a singular nonalgebraic point of a solution $w=w(z)$ of the above equation \eqref{eq:Painleverational}. There are two options:
either (i) $z_0$ is a transcendental point of $w(z)$ or (ii) $z_0$ is an essential singularity of $w(z)$.

Let us consider first (i). Then $w(z)$ either (ia) has a finite limit which can either (ia1)  be a zero of the equation
$$
Q(w, z_0)=0
$$
or (ia2) is not a zero of that equation, or (ib) tends to infinity. In (ia2) $Q(w_0, z_0)\ne 0$ and according to the Cauchy Theorem, there is a unique integral $w_0(z)$ determined by the initial condition $(w_0, z_0)$ and it is analytic in the neigbourhood of $z_0$. Thus, $w(z)$ coincides with $w_0(z)$ and is analytic at $z_0$.

In (ia1) case  $Q(w_0, z_0)= 0$ and thus  $P(w_0, z_0)\ne 0$, since $z_0$ is not in $\mathcal M$. According to the above Lemma \ref{lemma:lemma} about a critical point, $z_0$ is a critical algebraic point of the solution $w(z)$ of the equation  \eqref{eq:Painleverational}.

The case (ib) can be treated similarly, by changing variables to $w_1=1/w$.

Now, let us consider (ii): we assume that $z_0$ is an essential singularity of the integral $w(z)$. Denote by $W_j$ the zeros of the polynomial $Q(w, z_0)=0$. Let $D$
denote the region of undefinitness of the integral $w(z)$ at $z_0$, as the set of values the integral attains when the argument approaches the essential singularity.
Let $A$ be the point of $D$ closest to the origin (which can be the origin itself). Let $\bar w$ be  an arbitrary point of $D$ and $d_j$ the distances from $\bar w$ to $W_j$.
There are finitely many $d_j$s, they are all strictly positive, thus their minimum $d$ is strictly positive. Take any $\rho$ such that $0<\rho<d$ and construct circles $C_j$ centered at $W_j$ with radius $\rho$. Also take $R$, $R>|A|$ and construct circle $\Gamma$ centered at the origin with radius $R$. Let $D_0$ be the part of $D$ inside $\Gamma$ and outside all $C_j$.
Construct circles $C'_j$ concentric with $C_j$ with radius $\rho/3$. Select $r$ small enough that for all $z$ from the disk $S$ centered at $z_0$  with radius $r$, the solutions of he equation $Q(w, z)=0$ are inside $C'_j$. Let us also construct circles $C''_j$ concentric with $C_j$ with radius $\rho/2$ and circles $\Gamma'$, $\Gamma''$ concentric with $\Gamma$ with radii $R+\rho$ and $R+\rho/2$ respectively. Denote by $D_1$ the part of $D$ bounded by $\Gamma''$ and all $C_j''$.  There exists $M>0$ such that
\begin{equation}\label{ineq:estimate}
  \left|\frac{P(w, z)}{Q(w,z)}\right|<M
  \end{equation}
  on $S\times D_1$.
  Thus, for all $w_1$ in $D_0$ and $z_1$ in the disk $S_1$ concentric with $S$ and radius $r/2$, the disk $C$ of radius $\rho/2$ and center $w_1$ and the disk $\gamma$ with radius $r/2$ and center $z_1$ have the property that the function $P(w, z)/Q(w,z)$ is holomorphic on $\gamma\times C$ and satisfies the inequality \eqref{ineq:estimate}.

  From the Cauchy's result stated above, we know that the integral of the equation \eqref{eq:Painleverational}, defined by the initial condition $(w_1, z_1)$ is holomorphic within the disk centered at $z_1$ with radius
  $$
  \lambda=\frac{r}{2}(1-\exp(-\rho/(2Mr)))
  $$
  according to \eqref{eq:Cauchyestimate}.

  Now we can conclude the proof of the Painlev\'e Theorem in the case of the equations of the form \eqref{eq:Painleverational}.
  Indeed, we can construct a disk $\sigma$, centered at $z_0$ with radius $\lambda$. Then, within $\sigma$ there exists a point $z_1$ such that $w(z_1)=w_1$ belongs to $D_0$.
  Using the facts derived above, we see that the integral, analytic at $z_1$ and defined with $(w_1, z_1)$ as the initial conditions, is also holomorphic at $z_0$, since
  $$
  |z_1-z_0|<\lambda.
  $$
  Obtained contradiction concludes the proof.

Let us observe that, in the process of proving of  Theorem \ref{th:Painleve}, Painlev\'e also verified that the conditions of the L. Fuchs Theorem are also sufficient, as stated above, see Theorem \ref{th:Fuchs}. This remark about the L. Fuchs Theorem and the related Poincar\'e result is contained in the Section 8 of Chapter II of the First Part of the dissertation, on p. 57, see \cite{Painleve}, \cite{Painleve0}.

Paul Painlev\'e was born in Paris in 1863. He graduated from the \'Ecole Normale in 1877. He was a full professor of the { \'Ecole} Normale and Sorbonne. He was an elected member of the French Academy since 1900. After 1910, and election to the national parliament, Painlev\'e shifted his focus from science to politics. He was a minister
of several French governments, including the post of the Minster of War during the World War I. Painlev\'e served as  the Prime Minister of France two times: September 12 -- November 16, 1917 and April 17 -- November 28 1925. Painlev\'e died in Paris in 1933.

Let us conclude this Section by mentioning two important papers from 1889, the year when Petrovi\'c came to Paris, \cite{Kow} and \cite{Pic}.
The Kowalevski paper \cite{Kow} has appeared to become one of the most celebrated papers in the history of mathematics. Kowalevski successfully developed further some of the above ideas and concepts and applied them to the study of a system of algebraic equations, so called Euler-Poisson equations of motion of a heavy rigid body which rotates around a fixed point. Kowalevski investigated the possibility of such a system to have a general solution with poles as only possible movable singularities. In other words, she was looking for the cases where the general solution as a function of complex time is single-valued. As a result, she found what became to be known as the Kowalevski top. She integrated the equations of motion in that case explicitly using genus 2 theta functions and proved that her case indeed satisfied the initial analytic assumption. For these results Kowalevski received the famous Prix Bordin of the French Academy of Science, which was augmented from 3000 franks to 5000 franks. For the second paper on the rotation of a rigid body \cite{Kow1}, Kowalevski got a prize of the Swedish Royal Academy. Together with the later work of Painlev\'e and his students on the second order equations (see Section \ref{sec:conclusion}), these ideas of Kowalevski laid the  foundations of the so-called Kowalevski-Painlev\'e analysis, which is also known as a test of integrability. With this we conclude a brief description of the scientific atmosphere in which Petrovi\'c started the work on his PhD thesis.

\section{Petrovi\'c polygons}

The key preassumption of the main Petrovi\'c's construction is that a given point $x_0$ is a nonsingular point of the equation \eqref{eq1}.
In such a case, to each term  in the sum \eqref{eq2} with the coefficient $\varphi_i(x_0)={\rm const}\neq 0$ one corresponds a point in the $MON$ plane, according the following formulae
$$Q_i=(M_i,\,N_i), \quad M_i=m_{0i}+\dots+m_{ni}, \quad N_i=m_{1i}+2m_{2i}+\dots+nm_{ni}.$$
It should be noticed that one and the same point in the plane can correspond to one or more terms in the sum \eqref{eq2}. Let
$$S=\{Q_i,i=1\dots s,\;s\leqslant j\}$$
be the set of all points obtained in such a correspondence. We can draw these points in the $MON$ plane. In the Petrovi\'c dissertation  the set $S$ was extended with two segments $T_l$ and $T_r$  which are orthogonal to the axis $OM$ and connect the leftmost and the rightmost point of the set $S$ respectively with their projections to the $OM$ axis. The boundary of the convex hull of the set $S\bigcup T_l\bigcup T_r$ is a polygon. Both that polygon
and the concave part of the boundary of the convex hull of the set $S$ will be called  \textit{the Petrovi\'c polygon.} We will denote the Petrovi\'c polygon as $\mathcal{N}$. Let us point out that neither the vertical sides nor the horizontal side which lies at the $OM$-axis played any role in the applications of Petrovi\'c method. They don't correspond to any solution of the equation and in that sense their inclusion can be treated as artificial and unnecessary. However, they were included in Petrovi\'c's  original definition {\it not only} by pure formal or aesthetic reasons. They are indeed needed for methodological reasons as well, to allow a precise and elegant derivation of the main properties of the Petrovi\'c polygons. Thus, in our considerations, we will at the beginning, use the polygon as Petrovi\'c did, but later on in applications and computations, for simplicity we will not add these vertical segments any more and we will operate with the concave part of the convex hull of the set $S$ only. In our further deliberations, the Petrovi\'c polygon is an {\it irregular zig-zag line}. It is important to stress that, nevertheless, {\bf the conclusions from the considerations of the zig-zag line are identical to those coming from the entire polygon}. (Let us recall that Newton himself used irregular zig-zag lines, not polygons, \cite{Newton}, see also \cite{BK1986}, \cite{Ghys}.)

Let the equation \eqref{eq1} have a solution $y=y(x)$ which, in a  neighborhood of a point $x=x_0$, where $x_0$ is an arbitrary constant distinct from the singular points of the equation, can be presented in the form of  a power series with fractional exponents, i.e. in a form of the Puiseux series:
 \begin{equation}\label{eq3}
   y=\sum\limits_{k=0}^{\infty}c_k (x-x_0)^{\lambda+ k\sigma},
 \end{equation}
$\lambda=l\sigma,$ $\sigma\in\mathbb{Q}$, $\sigma>0$, $l\in\mathbb{Z}$, $c_k\in\mathbb{C}$, $c_0\neq 0$.
The main idea of Petrovi\'c was to use the polygon  $\mathcal{N}$ to keep those terms { $\varphi_i(x_0)\,y^{m_{0i}}y'^{m_{1i}}\ldots y^{(n) m_{ni}}$} of the equation \eqref{eq1} (called {\it the leading terms})
which form an equation (called {\it approximate equation}) having $\chi=c_0(x-x_0)^{\lambda}$  as its solution. Thus $\chi$ would be the asymptotic of the solution \eqref{eq3} in a  neighborhood of the point $x=x_0$. In that way we would effectively find such an asymptotic. As we see, these ideas of Petrovi\'c resemble the main idea behind the  Newton - Puiseux polygons in finding the asymptotics of solutions of algebraic equations. We are going to describe the ideas of Petrovi\'c in details, paying attention to the specifics of the case of ODEs.

Let us plug the formal series \eqref{eq3} into the equation \eqref{eq1}.
The solution  \eqref{eq3} and { its} derivatives can be rewritten in the form
 $$\begin{array}{lcl}
 y&=&(c_0+o(1))\;(x-x_0)^{\lambda},\\
 y'&=&(c_0\;\lambda+o(1))\; (x-x_0)^{\lambda-1},\\
 \hdotsfor{3}\\
 y^{(n)}&=&(c_0\;\lambda\cdot{\dots}\cdot(\lambda-n+1)+o(1))\; (x-x_0)^{\lambda-n}.
 \end{array}$$
Their powers have the following form:
$$\begin{array}{lcl}
 y^{m_{0i}}&=&\left(c_0^{m_{0i}}+o(1)\right)\;(x-x_0)^{\lambda m_{0i}},\\
 y'^{m_{1i}}&=&\left(c_0^{m_{1i}}\;\lambda^{m_{1i}}+o(1)\right)\;(x-x_0)^{\lambda m_{1i} -m_{1i}},\\
 \hdotsfor{3}\\
 y^{(n)m_{ni}}&=&\left(c_0^{m_{ni}}\;\left(\lambda\cdot{\dots}\cdot(\lambda-n+1)\right)^{m_{ni}}+o(1)\right)\;(x-x_0)^{\lambda m_{ni} -n m_{ni}}.
 \end{array}$$
 For every term of  $\varphi_i(x)\,y^{m_{0i}}y'^{m_{1i}}\ldots y^{(n)m_{ni}}$ for the formal series \eqref{eq3} we get the following expressions:
$$ \begin{array}{ll}\varphi_i(x)\,y^{m_{0i}}y'^{m_{1i}}\ldots y^{(n)m_{ni}}&=\varphi_i(x)\,(A_i(\lambda)\, c_0^{M_i} +o(1))\,(x-x_0)^{\lambda M_i-N_i} =\\
\\
&=(A_i(\lambda)\,c_0^{M_i}\,{\varphi_i}(x_0)+o(1))\,(x-x_0)^{\langle  R,\, Q_i\rangle},\end{array}$$ where
\begin{equation}\label{eqA}
A_{i}(\lambda)=\lambda^{\gamma_{1i}}(\lambda -1)^{\gamma_{2i}}\dots (\lambda -n+1)^{\gamma_{ni}},
\end{equation}
$$
\begin{aligned}
\gamma_{1i}={ m_{1i}+m_{2i}+\dots+ m_{ni},}\\
\gamma_{2i}=m_{2i}+m_{3i}+\dots+m_{ni},\\
\dots \dots \dots \dots \dots \dots \dots \dots \dots\\
\gamma_{ni}=m_{ni}.
\end{aligned}
$$
We will use the vector $ R=(\lambda,\,-1)$ and denote as $\langle\phantom{a}\rangle$ the dot product.

Thus, by substituting the series \eqref{eq3} into \eqref{eq1} we get the formula
 \begin{equation}\label{eq3_n1}
 (C_0+o(1))(x-x_0)^{\gamma_\lambda}=0,
\end{equation}
where the coefficient is given by
$$C_0=\sum\limits_{i:\,\langle  R, Q_i\rangle=\gamma_\lambda} A_i(\lambda)\, c_0^{M_i}\varphi_i(x_0),$$
where
$$\gamma_{\lambda}=\min_{i=1,\dots, s}\langle R,Q_i \rangle.$$
By a further analysis, formula \eqref{eq3_n1} can be rewritten in a more precise form as

\begin{equation}
\label{eq:formalsubstitut}
{C}_{k}(c_0, \dots, c_k)(x-x_0)^{\gamma_{\lambda}+k\sigma}{=0}.
\end{equation}
 Here $ C_{k}$ are polynomials of their arguments. Since the series \eqref{eq3} satisfies the equation \eqref{eq1}, the equation \eqref{eq:formalsubstitut} is satisfied identically, meaning that the coefficients $C_k$ in \eqref{eq:formalsubstitut} should all be zero. Consequently, if the solution of the equation   \eqref{eq1} exists in the form \eqref{eq3}, then by solving the equations $C_k(c_0, \dots, c_k) =0, k\in \mathbb Z_+$ we are getting the coefficients $c_k$ for which the series \eqref{eq3} gives the solution \eqref{eq1}. The {\it leading terms} of the equation    \eqref{eq1} are {those }
  $$ \varphi_i(x_0)\, y^{m_{0i}}\cdots y^{(n)m_{ni}}$$
 {corresponding} to the points $Q_i$, for which the dot product is {\it minimal}:
  $$\displaystyle \langle { R},\;Q_i \rangle=\langle (\lambda ,-1),\;(M_i,N_i) \rangle=\lambda\, M_i-N_i=\gamma_{\lambda}.$$

{
\bigskip
At the beginning of Chapter 1 of Part 1 of his thesis, Petrovi\'c proved the following statement.

 \bigskip
\begin{proposition}\label{prop:Petrovic} If $x=x_0$ is a nonsingular point of the equation { \eqref{eq1}}, and if $y=y(x)$ is a solution of the equation with an expansion  \eqref{eq3} with initial conditions  $y(x_0)=0$ or $y(x_0)=\infty$, then the first term  $c_0 (x-x_0)^{ \lambda}$ of the expansion into a Puiseux series of the solution \eqref{eq3} is a solution of the approximate equation, which corresponds either to a vertex or to a slanted edge of the polygon  $\mathcal{N}$ of the equation {\eqref{eq1}}.
\end{proposition}}

\begin{figure}[h]

    \qquad{\includegraphics[width=10cm]{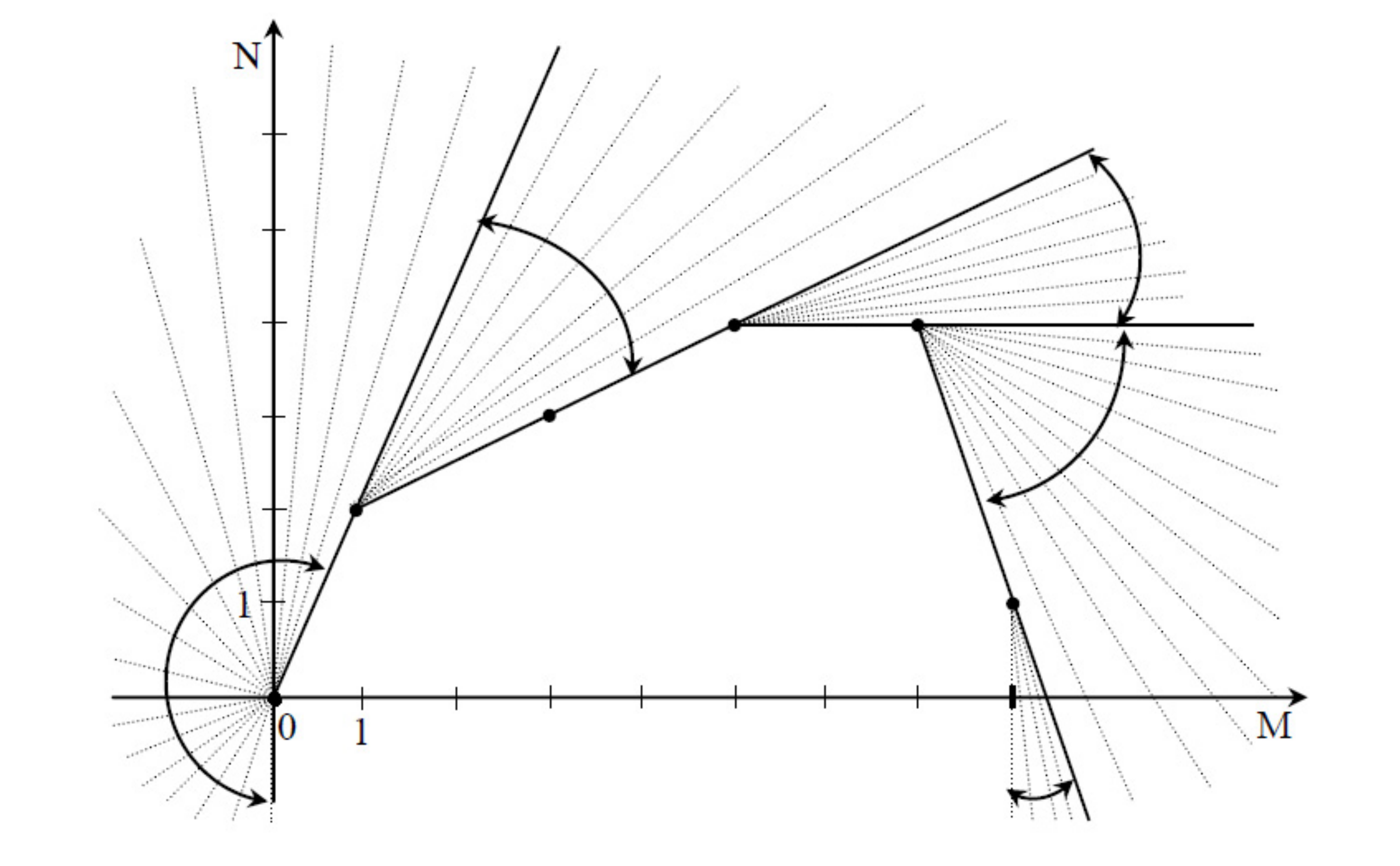}}
    \caption{Construction  of Petrovi\'c polygon}\label{fig:fig1}
\end{figure}

    Consider the line $\lambda M-N=\alpha$, $\alpha\in\mathbb{R}$, in the $MON$ plane. When the number $\alpha$ increases, the line $\lambda M-N=\alpha$ moves along the vector ${ R}$, directed inside the polygon. Thus, for the points $Q_i$ inside the polygon  $\mathcal{N}$  there is the relation for the dot product
    $$\langle { R},Q_i \rangle>\gamma_{\lambda},$$
     and its minimum attains at the boundary of the polygon: for {some} points $Q_i$ lying at the boundary of the polygon there is the following relation for the dot product
     $$\langle { R},Q_i \rangle=\gamma_{\lambda}.$$
      A point $Q_i$, such that $\langle { R},Q_i \rangle=\gamma_{\lambda}$, either is a vertex of the polygon or lies on one of its edges. Thus, the leading term  of the equation \eqref{eq1} corresponds either to a point which is a vertex of the polygon or to the points lying on an edge of the polygon.  It follows from the above constructions that there are finitely many values of $\lambda$, for which the minimum of the dot product $\gamma_{\lambda}$ attains on the edges, while there are infinitely many (continuum many) values of $\lambda$, for which the minimum of $\gamma_{\lambda}$ attains on vertices. The upper half-plane of the plane $MON$ can be decomposed on rays with the angular coefficients $\lambda$ and the open angular sectors, containing rays with angular coefficients $\lambda$, where $\lambda$ are all the values between the angular coefficients of the edges meeting at the given vertex, see Figure \ref{fig:fig1}.

Obviously, for any value of $\lambda$ there is a ray with the angular coefficient $\lambda$, which either corresponds to an edge or to a vertex. Thus, if the equation  \eqref{eq1} has a solution $y(x)$, represented in a neighborhood of the point $x=x_0$ in the form of the series \eqref{eq3},  then this solution has to correspond either to a vertex or to an edge of the polygon. Moreover, if $\lambda<0$, which means that  $x=x_0$ is a pole of $y(x)$, then these values correspond either to edges or vertices of {\it the right part} of the polygon $\mathcal{N}$. Similarly, if $\lambda>0$, i.e. $x=x_0$ is a zero of the solution $y(x)$, then corresponding edges and vertices belong to {\it the left part} of the polygon  $\mathcal{N}$.

In our clarifications of the construction of the Petrovi\'c polygon, we have used the Puiseux series. If, instead of the Puiseux series one considers generalized power series, i.e. the series with complex exponents (the complex exponents belong to a finitely generated semi-group), then the idea and the mode of the construction of the polygon remains the same.  Let us notice that Petrovi\'c in his thesis considered complex exponents of the asymptotic of the solutions in a neighborhood of nonsingular points of an algebraic ODE.

Let us clarify how to detect the complex exponents by using the Petrovi\'c polygon. First observe, that since the given equation is algebraic, the edges have rational slopes, i.e.
 the minimum of the dot product on edges is attained only for { the rational  exponents}, while on vertices one can have { complex exponents of asymptotics}. Let us consider now the general case, when one vertex  $Q_1=(\alpha,\beta)$ of the polygon $\mathcal{N}$ corresponds to more than one term in the equation \eqref{eq1}, i.e. in the sum
\begin{equation}\label{eq16}
\sum\limits_i\varphi_i(x)\; y^{m_{0i}}y^{\prime m_{1i}}y^{\prime \prime m_{2i}}\dots y^{({ n})m_{{ n}i}},
\end{equation}
$$
\varphi_i(x_0)={\rm const}\neq 0, \, \alpha=m_{0i}+m_{1i}+m_{2i}+\dots+m_{{ n}i}, \, \beta=m_{1i}+2m_{2i}+\dots+{ n}m_{{ n}i}.
$$
We expand the sum \eqref{eq16} in series in powers of $x-x_0$ and then into that we substitute the { expression} $y={ (c_0+o(1))(x-x_0)^{\lambda}}$ (here $c_0\neq 0$ { is an arbitrary constant} and $\lambda$ { is a sought exponent}), to get the expression
   $$\left(c_0^{\alpha} \sum_i \varphi_i(x_0)\;\lambda^{m_{1i}}\cdots(\lambda(\lambda-1)\cdots(\lambda- n+1))^{m_{ni}}+o(1)\right)(x-x_0)^{\lambda\alpha-\beta}.$$
    The leading terms in this sum are the following ones:
    \begin{equation} c_0^{\alpha} \sum\limits_i\varphi_i(x_0)\;\lambda^{m_{1i}}\cdots(\lambda(\lambda-1)\cdots(\lambda- n+1))^{m_{ni}}\;(x-x_0)^{\lambda \alpha-\beta}.\label{eq10}\end{equation}
     {As in \eqref{eqA}}, we introduce the following notation
$$
A_{i}(\lambda)=\lambda^{\gamma_{1i}}(\lambda -1)^{\gamma_{2i}}\cdots (\lambda -n+1)^{\gamma_{n i}},
$$
where
$$
\begin{aligned}
\gamma_{1i}=m_{1i}+m_{2i}+\dots+ m_{ni},\\
\gamma_{2i}=m_{2i}+m_{3i}+\dots+m_{ni},\\
\dots \dots \dots \dots \dots \dots \dots \dots \dots\\
\gamma_{ni}=m_{ni}{,}
\end{aligned}
$$
{and rewrite} the expression \eqref{eq10} in the form
 $${c_0^{\alpha}}(x-x_0)^{\lambda \alpha-\beta}\sum\limits_{i}A_i({\lambda})\varphi_i(x_0).$$
 This sum is equal to zero only if the polynomial
 $$\sum\limits_{i}A_i({\lambda})\;\varphi_i(x_0)=a_{\beta} \lambda^{\beta}+\dots+a_1\lambda+a_0$$
 is zero.

 \begin{definition}[Petrovi\'c, 1894] The equation
 \begin{equation}\label{eq12}
 a_{\beta} \lambda^{\beta}+\dots+a_1\lambda+a_0=0
\end{equation}
 will be called {\it the characteristic equation} of a given vertex.
 \end{definition}

 Obviously, if  $\lambda\in\mathbb{C}$ satisfies the characteristic equation  \eqref{eq12} and if it is the exponent of the first term of the expansion of a solution of the given ODE, then the minimum of the dot products $\min\limits_i \langle({\rm Re\,}\lambda,-1),Q_i\rangle=\gamma_{\lambda}$ attains only on the vertex  $Q_1$, i.~e.
 $$\langle({\rm Re\,}\lambda,-1),Q_1\rangle=\gamma_{\lambda}.$$

For an algebraic ODE of the first order, according to the Painlev\'e Theorem \ref{th:Painleve}, its solutions could possess only algebraic movable singular points, i.e.  every solution
$y=y(x)$ in a neighborhood of a nonsingular point $x=x_0$ of the equation presents in a form of a power series \eqref{eq3} in a general case with fractional exponents and uniquely defined coefficients. By using a generalized method of Newton-Puiseux polygons, the method of polygons of Petrovi\'c or of Fine, one can completely determine the expansions of the solutions of an algebraic ODE of the first order around a point which is nonsingular for the equation.  Because of that, Petrovi\'c was able to fully resolve the question about the conditions under which the solutions of an algebraic ODE of the first order have or not have movable zeros or movable poles, see Theorems \ref{th:1}-\ref{th:3}.

We should make an important observation. In his studies of algebraic ODEs of the first order, Petrovi\'c considered only slanted  edges of a polygon.
If a considered solution has the initial data $y(x_0)=0$ or $y(x_0)=\infty$, then, obviously, it corresponds to slanted edges of the Petrovi\'c polygon, left or right { respectively.
The vertices and horizontal edges} of the Petrovi\'c polygon correspond to the solutions of an algebraic ODE with the initial data
  $y(x_0)=C={\rm const}\neq 0$. In this case, the point  $x=x_0$ is not a zero. By making a change of dependent variable  $y=C+u$ in the initial equation, and then, by using the method of Petrovi\'c polygon to the transformed equation, one can determine the order of the zero of the solution $u=u(x)$, $u(x_0)=0$.

Let us notice that if $x=x_0\,$ is not  a  singular point of a given algebraic ODE of the first order, then it can be:
  \begin{itemize}
  \item
  either  a nonsingular point of the solutions, and in that case there is a one-parametric in $C$ family of holomorphic solutions in a neighborhood of that point; in other words { the assumptions of the implicit function theorem and Cauchy theorem}  are met;
     \item
      { or a movable  algebraic singular point} of those solutions; in this case all the coefficients of the corresponding expansions { in  Puiseux}  series with a finite principal part,  are uniquely determined.
     \end{itemize}

     \begin{example}\label{ex:ex1} Let us consider the equation
  \begin{equation}\label{eq14}
  f(y',y)=y'^2(y-1)+1=0.
\end{equation}
It possesses two general solutions
$$y=\left(\frac{C\pm 3x}{2}\right)^{2/3}+1.$$

Without calculating the solutions of the equation \eqref{eq14}, one can observe that the solutions have no {fixed} singular points. This is a consequence of the fact that the coefficients of the equation have no { zeros}.  { Moreover, for an arbitrary  $x_0\in\mathbb{C}$, solutions with the initial condition $y(x_0)\neq 1$ are holomorphic  in a neighbourhood of the point $x=x_0$, and solutions with the initial condition  $y(x_0)=1$ (which correspond to integral constants $C=\pm 3x_0$) have $x_0$ as a movable algebraic critical singular point.}

 Let us investigate the equation \eqref{eq14} by use of the Petrovi\'c polygon. The terms $y'^2y$, $-y'^2$, $1$ of the equation \eqref{eq14} correspond to the points  $(3,2)$, $(2,2)$, $(0,0)$.  The polygon of the equation \eqref{eq14}  (see Figure \ref{fig:fig2}a) has two edges: the horizontal one  $[(2,2),\;(3,2)]$ and a left-sloped one $[(2,2),\;(0,0)]$ with the angular coefficient equal to 1. There are three vertices: $(3,2)$, $(2,2)$, $(0,0)$.

\begin{figure}[h]

    \qquad{\includegraphics[width=12cm]{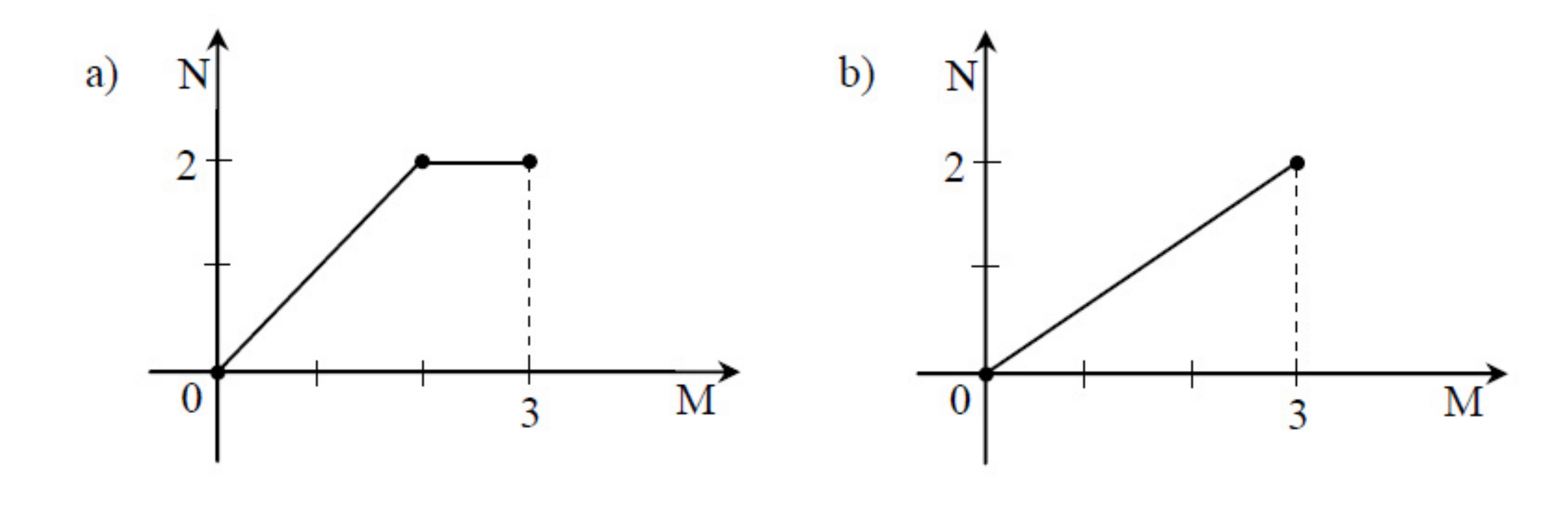}}
    \caption{The Petrovi\'c polygons of: a) the equation \eqref{eq14}; b) the equation \eqref{eq14b}.}\label{fig:fig2}
\end{figure}

 The vertices $(3,2)$ and $(2,2)$ correspond to approximate equations  $y'^2y=0$ and  $y'^2=0$ respectively. The approximate equations correspond to the solutions  $y={\rm const}\neq 0$, i. e. the exponent is  $\lambda =0$. When we substitute into the equation an expansion of the solution, which begins with a constant, we observe that the minimum of the dot products of the exponents
 $$\min\limits_{i} \langle(\lambda,-1),Q_i\rangle=\min\left( \langle(0,-1),(3,2)\rangle,\; \langle(0,-1),(2,2)\rangle,\; \langle(0,-1),(0,0)\rangle\right)=-2$$
  attains not at a single vertex, $(3,2)$ or $(2,2)$ but along the entire horizontal edge  $[(2,2),\;(3,2)]$, i.e. $\langle(0,-1),(3,2)\rangle=-2$ and $\langle(0,-1),(2,2)\rangle=-2$. Thus, there are no solutions here which would correspond to the vertices.
The solution which begins with a constant corresponds to the horizontal edge $[(2,2),\;(3,2)]$ with the approximate equation  $y'^2(y-1)=0$, which possesses two solutions $y=C_0={\rm const}\neq 0$ and $y=1$. By changing the dependent variable $y=C_0+u$, { $C_0\neq 1$,} in the equation \eqref{eq14} we get the equation $u'^2(u+C_0-1)+1=0$ with the same polygon as the equation \eqref{eq14}. However, there is now an additional condition  $u(x_0)=0$. We should not consider the vertices of the polygon in this new situation because the approximate equations which correspond to them lead to constant solutions, in other words $u(x_0)\neq 0$. Since the expansion of a solution is given in increasing order of powers of  $x-x_0$, we are interested in edges with positive slopes only, these are the left-sloped edges. In the new polygon, this is the edge with the slope 1. { This edge } corresponds to the approximate equation $u'^2(C_0-1)+1=0$ with the solution $ u=\pm (x-x_0)/\sqrt{-C_0+1}$. Obviously, { according to the Cauchy theorem the solution}
\begin{equation}y=C_0\pm (x-x_0)/\sqrt{-C_0+1}+\dots,\qquad C_0\neq 1,\label{eq15}\end{equation}
of the equation \eqref{eq14} {  coincides with the holomorphic solution that satisfies the initial condition $y(x_0)=C_0$.} If we now change the dependent variable $y=1+u$, $u(x_0)=0$
in the equation  \eqref{eq14}, we come to the equation
\begin{equation}\label{eq14b}
u'^2u+1=0.
\end{equation}
 The polygon of this equation consists of only one, slanted edge $[(3,2),\;(0,0)]$ with the positive slope equal to $2/3$, and two vertices  $(3,2)$ and $(0,0)$ (see Figure \ref{fig:fig2}b). Again, we don't consider the vertices of the new polygon, since they correspond only to solutions of the new equation with expansions with  $u(x_0)\neq 0$. The edge $[(3,2),\;(0,0)]$ corresponds to the approximate equation $u'^2u+1=0$ with the solution $u=\pm(3/2(x-x_0))^{2/3}.$ Obviously, with $x_0$ fixed, we get uniquely determined  solution
$$y=1\pm(3/2(x-x_0))^{2/3},$$
which is a particular solution of the equation \eqref{eq14} { with the initial condition $y(x_0)=1$}. { Here the point  $x_0$ is a movable critical algebraic singular point of the solutions. But, it is not a zero of those solutions.}

There is one more element of the polygon of the equation \eqref{eq14}, which remains to be considered: the slanted edge $[(2,2),\;(0,0)]$  with the angular coefficient $1$. This edge corresponds to  the approximate equation  $-y'^2+1=0$ with the solution $y=\pm (x-x_0)$. Clearly, the solution  $y=\pm (x-x_0)+\dots$ of the equation  \eqref{eq14} is holomorphic according to { the Cauchy theorem and it is a particular solution of the family of solutions \eqref{eq15} for $C_0=0$.} We conclude this simple example which illustrates the method of investigation of the singularities of the solutions of an algebraic ODE based on the use of the Petrovi\'c polygon.
\end{example}

\subsection{Generalized homogeneity and  some limitations  concerning higher order ODEs}\label{sec:obstructions}

\bigskip
 Let us recall that an approximate equation is called {\it generalized homogeneous} in $x-x_0$ (or in  $y$), if it is invariant with respect to the change of $x-x_0$ to $k(x-x_0)$  (of  $y$ to $ky$), $k\in\mathbb{C}$. Approximate equations which are generalized homogenous often can be solved explicitly.
An important feature of the approximate equations, obtained through the Petrovi\'c polygons is their {\it generalized homogeneity}. This means that if there exists an approximate equation
 $$h(x_0,y,\dots,y^{(n)})=0,$$
  corresponding to a slanted edge with the angular coefficient equal to $\lambda$, then after the transformation  $y=(x-x_0)^{\lambda}u$ the approximate equation transforms to a new equation $(x-x_0)^{\lambda}\widetilde{h}(x_0,u,(x-x_0)u',\dots, (x-x_0)^nu^{(n)})=0$, where  $\widetilde{h}$ is a polynomial function of its arguments and a generalized homogeneous
  function in  $x-x_0$.

  If an edge is horizontal, then the corresponding approximate equation defines  a generalized homogeneous function in $x-x_0$. If an edge is vertical, the corresponding equation is a generalized homogeneous in $y$.  An approximate {equation}  which corresponds to a vertex is generalized homogeneous both in $x-x_0$ and in  $y$.

As a benefit from the Painlev\'e Theorem \ref{th:Painleve}, Petrovi\'c did not need to consider a higher-dimensional constructions -- the polyhedra of the algebraic ODE of the first order in order to investigate fully the movable singularities of their solutions. According to the Painlev\'e Theorem all movable singularities of such equations are algebraic and the planar polygon captures all such singularities.   It is very important to stress that Petrovi\'c in his dissertation clearly pointed out the limitations of the applications of his polygonal method to the algebraic ODEs of higher orders. He showed that the method of planar polygons could be successfully applied to higher order algebraic ODEs to study some types of movable singularities of the solutions. However, due to the lack of a Painlev\'e type result for higher order equations, Petrovi\'c understood that his method was powerless in proving absence of other types of movable
singularities. In other words, the algebraic ODEs starting from order two can have movable singularities which are not algebraic. Moreover, the algebraic ODEs starting from order three can even have non-isolated movable singularities. If we pass from the Petrovi\'c polygons to higher dimensional polyhedra with the aim to study { non-algebraic} movable singular points of algebraic ODEs of higher order, we may at first hope to use the results of the modern theory of Newton polyhedra. However, the approximate equations obtained through the polyhedra are quite complex, they don't possess the property of generalized homogeneity, and often are not exactly solvable.

\bigskip

\section{Fine polygons}

Fine as well generalized the polygonal method of Newton and Puiseux. He used his generalization to study formal asymptotics of the solutions of algebraic ODEs
\eqref{eq1} at the point $x=0$. In his considerations, he includes both cases, when the point $x=0$ is a singular point of the equation and also when it is not a  singular point
of the equation. In his papers \cite{Fine}, \cite{Fine1}, Fine  used Puiseux \cite{Puiseux} and Briot-Bouquet \cite{BriotBouquet1} results and generalizes them. Fine and Petrovi\'c methods of construction of approximate equations are based on the same principles. Therefore, this is natural that the Fine method matches the steps of the Petrovi\'c method.
 In the construction of Fine polygons, we correspond a point to every term of the equation of the type
  $$ c\,{x^{l_{it}}}y^{m_{0i}}y'^{m_{1i}}\cdots y^{(n) m_{ni}}, \qquad c\in\mathbb{C},$$
  where the point  $({N}_{it},M_i)$, is determined by the formula
${N}_{it}=l_{it}-N_i,$ where  $N_i$ and $M_i$ are defined in the same way as in Petrovi\'c's construction above. If the points  $({N}_{it},M_i)$ are depicted in the plane and if we consider the boundary of the convex hull of all the points  $({N}_{it},M_i)$, then {\bf the left part of that boundary} (consisting of the edges and vertices where the external normal is pointed to the left) captures the behavior of the solutions in a neighborhood of the point $x=0$. We will call this left part of the boundary {\it the Fine polygon}. The vertices and edges of the Fine polygon correspond to the leading terms of the equation, i.e. those terms of the equation \eqref{eq1} which can form approximate equations. The candidates for the role of the asymptotics of the true solutions of the original equation lie among the solutions of the approximate equations.  Let us observe that the Fine polygon takes into account the exponents  $l_{it}$ of the independent variable $x$ in the coefficients  $\varphi_i(x)$ of the equation \eqref{eq1}, because here $x=0$ can be a singular point for the equation \eqref{eq1}, i. e.  $\varphi_i(0)$  could be zero or be undefined.

Let us also observe that by using the change $x=z+x_0$, the problem of analysis of the solutions  in a neighborhood of an arbitrary point $x=x_0$ reduces to the problem of the analysis of solutions in a neighborhood of the point $z=0$. Thus, we proved the following:

\begin{theorem}
The Fine polygon of the equation $f(z+x_0,{y,\,}y',\dots,y^{(n)})=0$, where  $x_0$ {\bf is not a singular point of the equation} \eqref{eq1}, coincides with the Petrovi\'c polygon of the equation \eqref{eq1} rotated by $\pi/2$ in the counterclockwise direction.
\end{theorem}

Fine's paper \cite{Fine} is mostly devoted to the question of calculations of terms in the expansion of formal solutions, (which have a form of Puiseux series) of algebraic ODEs in a neighborhood of the point $x=0$. Fine also treated the question of the convergence of formal series. Fine proved the following result:

\begin{theorem}[Fine,  \cite{Fine}]
If every term of an algebraic ODE contains {the dependent variable and its derivatives of all orders}, i.e. if every term
$$
\varphi_i(x)\,y^{m_{0i}}y'^{m_{1i}}\cdots y^{(n)m_{ni}}
$$
of \eqref{eq1} and \eqref{eq2} satisfies $m_{0i}, \dots, m_{ni}>0$,
then all the formal Puiseux series which formally satisfy the given equation
converge.
\end{theorem}

\section{Further development of the methods of polygons of Petrovi\'c and Fine}

 A century after this Fine's result, Malgrange \cite{Malgrange} gave {sufficient} conditions for convergence of a formal power solution, which solves a given analytic ODE.

As it has been mentioned in Introduction,  J. Cano developed further the method of Fine polygons in \cite{Cano2}, \cite{Cano1}, \cite{Cano3}. He applied these methods in calculations of formal solutions of ODEs and partial differential equations and to prove the convergence of the formal solutions.

Not being aware of the works of other authors, the predecessors, Fine and Petrovi\'c, and his own contemporaries, Cano, Grigor'ev, and Singer (\cite{GrigorievZinger}), A. D. Bruno suggested the methods to calculate formal solutions of algebraic differential equations and systems of equations (see \cite{Bruno}, \cite{Bruno1}). These methods are also based on generalizations of the Newton-Puiseux polygons and they repeat the ideas of Petrovi\'c and Fine, enriching them with some additions and extensions. These additions and extensions are, essentially, related to calculations of finitely-generated semi-groups of exponents of terms of generalized formal series, which formally satisfy a given algebraic ODE and also to extensions of the classes of the considered formal solutions.

To explain these ideas, let us introduce the notion of the order  $p(f)$ of function $f(x)$. If there exists the limit
 $$
 \lim\limits_{x\rightarrow 0, x\in \mathcal D}\frac{\ln{|f(x)|}}{\ln{|x|}}=\lim\limits_{x\rightarrow \infty, x\in \mathcal D}\frac{\ln{|f(x)|}}{\ln{|x|}}=p(f, \mathcal D)\in\mathbb{R}\bigcup\{\pm \infty\},
 $$
  then the value $p(f, \mathcal D)$ will be called  {\it the order of the function $f(x)$ with respect to a domain $\mathcal D\subset \mathbb C$}, where the closure of $\mathcal D$ contains the points $0$ and  $\infty$.  Similar definition of the order of a function can be found in the work of Nevanlinna   \cite{Nevanlinna}. Notice that the power functions, logarithms, and the products of such functions have finite orders. Moreover, the order of the derivatives of these functions, if {the derivatives are} not identically equal to zero, with each differentiation decreases its order by 1. For example $p(x^2\ln\ln x)=2$, and $p((x^2\ln\ln x)')=1$, $x\in\mathbb{C}$. The same rule does not work for the functions {$\cos x$, $x\in\mathbb{C}$, $0<|{\rm Im}\,x|<A$,} since: $p(\cos x)=0$ è $p((\cos x)')=0$. The constructions of polygons of Petrovi\'c and Fine take into account the finiteness of the order and the rule of change of the order by 1 with each differentiation of a power function or of a formal series.   We will use the sign $\leftrightarrow$
  to denote the correspondence between a term of a given algebraic ODE and a point of its polygon in the plane. For simplicity, let us consider the terms  $y,\;y',\;\dots,\;y^{(n)}$. Indeed, Petrovi\'c had used the following correspondence:
$$y\leftrightarrow(1,0),\;y'\leftrightarrow (1,1),\dots,y^{(n)}\leftrightarrow(1,n),$$
while Fine's choice was:
$$y\leftrightarrow(0,1),\;y'\leftrightarrow (-1,1),\dots,y^{(n)}\leftrightarrow(-n,1).$$
 Taking into account various classes in which one searches for formal solutions of ODEs, it is possible to generalize the polygons of Petrovi\'c and Fine.  One of that was accomplished in  \cite{BrunoGoryuchkina2}. See also Section \ref{sec:conclusion} for some further considerations.

\section{On movable singularities of algebraic ODEs of the first order}

In the first part of the dissertation, Petrovi\'c considers algebraic ODEs of the first order of a general type:
\begin{equation}\label{eq4}
  F(x,y,y')=0,\qquad F=\sum\limits_{i=1}^{j} \varphi_i(x)y^{m_{0i}}y'^{m_{1i}}.
\end{equation}

\bigskip
In the sequel, the notions of zeros and poles will include both ordinary and critical zeros and poles.

 \bigskip
   Consider a point $x=x_0$, where $x_0$ is not a singular point of the equation  \eqref{eq4}. Assume that the equation \eqref{eq4} has a solution $y=y(x)$, such that $y(x_0)=0$ or $y(x_0)=\infty$. According to the Painlev\'e Theorem \ref{th:Painleve}, the expansion of the solution $y=y(x)$ in a  neighborhood of the point  $x=x_0$ is a power series, in a general case, with fractional exponents. All such solutions are detectable through the method of Petrovi\'c { polygon whose distinctive feature in the case of an algebraic ODE of the first order is that each its vertex corresponds to the exactly one monomial $\varphi_i(x_0)y^{m_{0i}}y'^{m_{1i}}$.}

\bigskip
 At the beginning of Chapter 1 of Part 1 of his thesis, after Proposition \ref{prop:Petrovic}, Petrovi\'c formulates and proves necessary and sufficient conditions for the absence of movable zeros and poles of the solutions of the equation \eqref{eq4}.

\bigskip
 \begin{theorem}[Petrovi\'c, 1894] \label{th:1} The necessary and sufficient condition for poles $($zeros$)$ of the general solution of a given algebraic ODE of the first order \eqref{eq4} not to depend on the constants of integration is that the polygon  $\mathcal{N}$ of the equation \eqref{eq4} does not contain right $($left$)$ slanted edges.
 \end{theorem}

\bigskip {
Sufficiency in Theorem \ref{th:1} follows from Proposition \ref{prop:Petrovic} whereas necessity  is proved by the methods of analytic theory of differential equations. In particular Petrovi\'c used some techniques of L. Fuchs from \cite{Fuchs} from the proof of what we presented as Theorem \ref{th:Fuchs}.}

 \bigskip
 Thus, the situation, when both zeros and poles of the general solution of the algebraic differential equation of the first order  \eqref{eq4} are not movable, can appear if the polygon $\mathcal{N}$ of the equation \eqref{eq4} is a horizontal segment or more generally  a part of { a rectangle.

 Let us underline that the conditions} of Petrovi\'c Theorem \ref{th:1} can be easily verified just through a simple construction of the polygon of the equation, while the verification of the conditions of the Theorem of L. Fuchs require much more elaborated work which includes the discriminant equation to be solved, then the equation   $F(x,y,y')=0$ to be resolved with respect to the derivative, and finally
an expansion of the derivative in a  series in a neighborhood of a discriminantal solution is needed.
However, Petrovi\'c's Theorem allows to promptly check the existence or the absence of movable zeros or poles, while a situation of the absence of movable critical singularities which are not zeros, could be detected either by the use of the Fuchs Theorem or by repeated use of Petrovi\'c's Theorem, see Example \ref{ex:ex1} and Examples \ref{ex:ex2}, \ref{ex:ex3} below.

\bigskip
Further on in Chapter 1 of Part 1 of the thesis, Petrovi\'c continues to study the existence of singular points of the general solution of the equation
 \eqref{eq4}. He formulates and proves a few theorems.

  \bigskip

\begin{theorem}\label{th:2}
The necessary and sufficient conditions for the general solution of a given algebraic ODE of the first order   \eqref{eq4} to possess a movable zero (or a pole) of order  $\lambda$ is that the polygon $\mathcal{N}$ of the equation \eqref{eq4} contains an edge with the angular coefficient $\lambda$ $($or $-\lambda)$.
\end{theorem}

\bigskip
\begin{theorem}\label{th:3}
If the general solution of a given algebraic ODE of the first order   \eqref{eq4}  possesses poles which are independent of the constants of integration, then there are finitely many such poles.
\end{theorem}

\bigskip

\begin{example}\label{ex:ex2} Consider the equation
\begin{equation}\label{eq5n}
xy'^3+yy'-1=0.
 \end{equation}
 This equation solves implicitly, for example with the assistance of computer algebra. Obviously there are three general solutions since the left hand side of the equation
 \eqref{eq5n} is a polynomial of degree three in $y'$. These solutions are cumbersome and because of that we are not going to list them here. Let us check the conditions of the Fuchs Theorem for the equation  \eqref{eq5n}. As we see, the first two conditions of Fuchs Theorem \ref{th:Fuchs} are satisfied. The discriminant equation  $y^3 +27 x/4=0$  has a multi-valued solution $y=\sqrt[3]{-27x/4}$, which is not a solution of the equation  \eqref{eq5n}. The third condition of the Fuchs Theorem is not satisfied. Thus, solutions of the equation \eqref{eq5n} possess movable critical points.

 Let us now turn to Petrovi\'c's Theorems \ref{th:1} and \ref{th:2}.
The points  $Q_1=(3,3),$ $Q_2=(2,1)$, and $Q_3=(0,0)$ correspond to the equation \eqref{eq5n}. The Petrovi\'c polygon  (see Figure \ref{fig:fig3}) consists of one edge $[(3,3),$ $(0,0)]$ with the angular coefficient  equal to $1$. According to Theorem \ref{th:2} this edge corresponds to the solutions of the equation  \eqref{eq5n} with movable zeros of order $1$. Solutions with movable zeros of order greater than $1$ do not exist.

\begin{figure}[h]

    \qquad{\includegraphics[width=10cm]{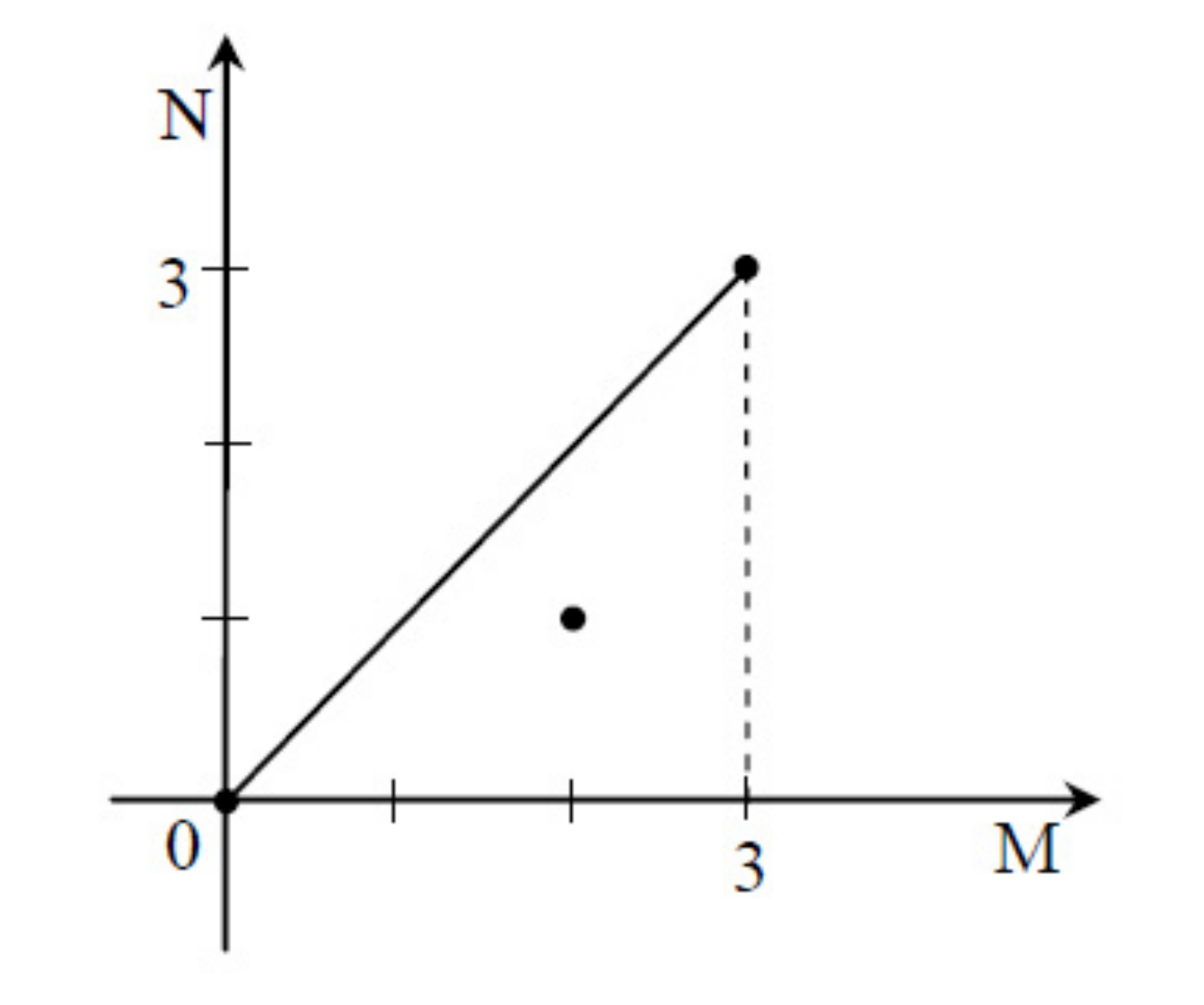}}
    \caption{The Petrovi\'c polygon of the equation \eqref{eq5n}.}\label{fig:fig3}
\end{figure}

 According to { the Cauchy theorem} in a  neighborhood of a nonsingular point
 $x=x_0$ ($x_0\neq 0,\infty$) of the equation \eqref{eq5n} any of its solutions having zero of order one at that point, $y(x_0)=0$, is presented by a uniformly convergent series
 \begin{equation}y=c_1 (x-x_0)+\sum\limits_{k=2}^\infty c_k (x-x_0)^k,\label{eq13n} \end{equation}
 where $c_1=\sqrt[3]{1/x_0}$, and the consecutive complex coefficients  $c_k$ are  uniquely determined  coefficients. Thus, all movable zeros of equation \eqref{eq5n} are noncritical, and the existence of movable critical points follows from the Fuchs Theorem.
\end{example}

\bigskip
\begin{example}\label{ex:ex3} Consider the equation
\begin{equation}\label{eq5}
f(x,y,y')=y'^2-(y'-1)(y-1)+x=0.
 \end{equation}
 It has two singular points $x=1$ and $x=\infty$. It solves implicitly.  It has two general solutions since  $f(x,y,y')$ is a polynomial of the second degree in $y'$.
 These solutions are cumbersome and because of that we are not going to list them here. Let us check the conditions of the Fuchs Theorem for the equation
 \eqref{eq5}. The first  condition of Fuchs Theorem \ref{th:Fuchs} is satisfied. The second condition is not satisfied.  Thus, solutions of the equation \eqref{eq5}  possess movable critical points.

   We want to use Petrovi\'c's Theorems \ref{th:1} and \ref{th:2}. The points  $Q_1=(2,2),$ $Q_2=(2,1)$, $Q_3=(1,1)$, $Q_4=(1,0)$, $Q_5=(0,0)$  correspond to the equation \eqref{eq5}. The Petrovi\'c polygon (see Figure \ref{fig:fig4}) consists of one slanted edge  $[(2,2),$ $(1,1),$ $(0,0)]$  with the angular coefficient  equal to $1$. According to Theorem \ref{th:2} this edge corresponds to the solutions of the equation  \eqref{eq5} with movable zeros of order $1$, and movable zeros with order higher than $1$ are absent.

\begin{figure}[h]

    \qquad{\includegraphics[width=10cm]{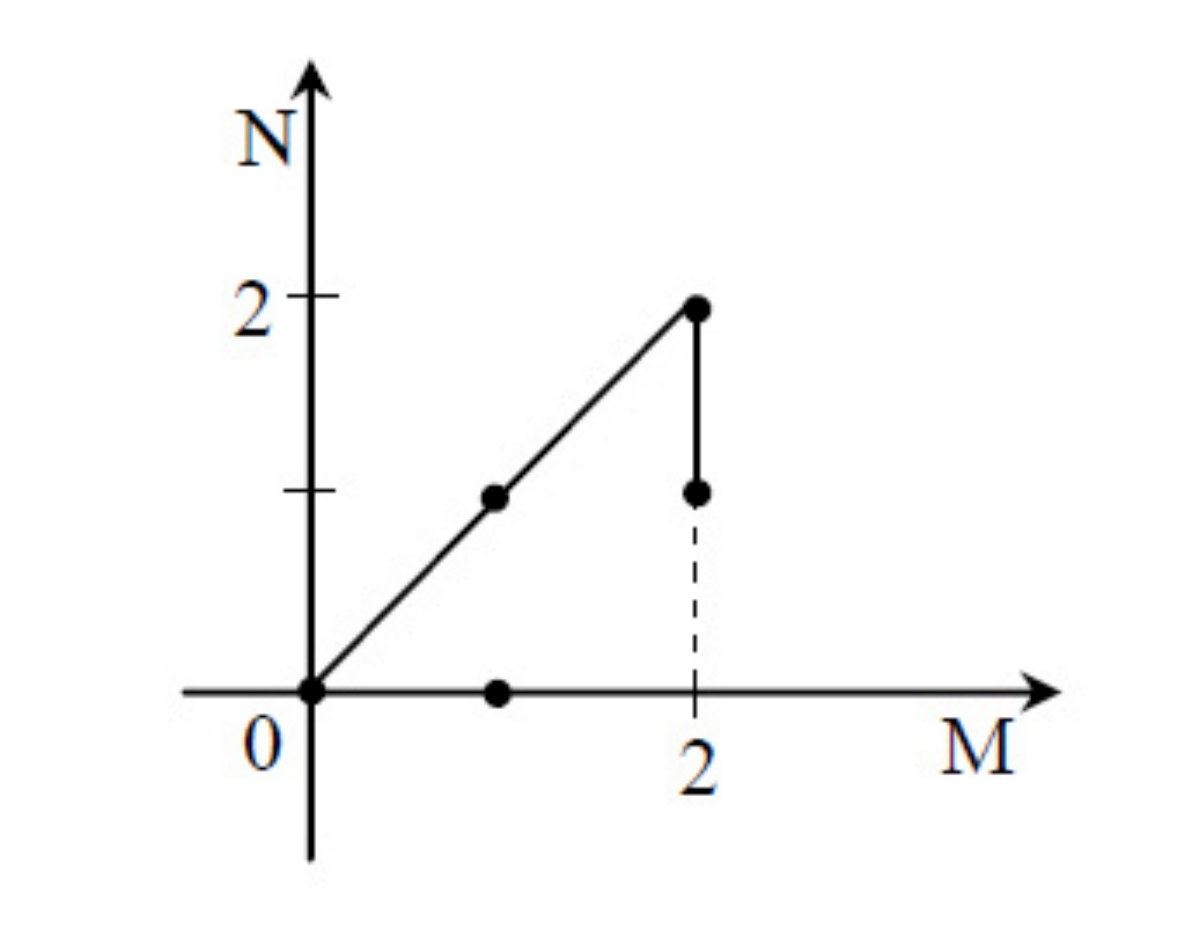}}
    \caption{The Petrovi\'c polygon of the equation \eqref{eq5}.}\label{fig:fig4}
\end{figure}

 Observe that if  $y(x_0)=0$, where  $x_0$ is a nonsingular point of the equation i.e. $x_0\neq 1,\;\infty$, then $y'(x_0)\neq 0$.

 In the case  {$x_0\neq 1,\,\frac54,\,\infty$, $y(x_0)=0$,} the Cauchy theorem can be applied to the equation  \eqref{eq5}, and thus, there exists a solution which is represented by a uniformly convergent series
 $$y=c_1 (x-x_0)+\sum\limits_{k=2}^\infty c_k (x-x_0)^k,$$
where
 $$\displaystyle c_1=-\frac{1\pm\sqrt{5-4x_0}}{2},$$
  and other complex coefficients  $c_k$ are  uniquely  determined coefficients.

   However, the Cauchy theorem is not applicable in the case $x_0=\frac54$, $y(\frac54)=0$, $y'(\frac54)=-\frac12$ since $f(\frac54,0,-\frac12)=0$, $\frac{\partial{f}}{\partial{y'}}(\frac54,0,-\frac12)=0$. In the neighborhood of the point $x_0=\frac54$ there exists a solution of the equation  \eqref{eq5},
   which can be presented by the Puiseux series
   $$y=-\frac{1}{2} \left(x-\frac54\right)+\sum\limits_{k=2}^\infty c_k \left(x-\frac54\right)^{(k+1)/2}.$$
   The Petrovi\'c polygon method is, in that sense, more universal than the Cauchy theorem. It is able to investigate all solutions with power asymptotics in a neighborhood of a point which is nonsingular for the equation. In this way, there is a movable critical zero of the solution in a neighborhood of a point ($x_0=5/4$) which is not singular for the equation.
\end{example}

\section{On single-valued solutions of algebraic ODEs of the first order explicitly resolved w.r.t. the derivative. Generalized Riccati equations}

Let us recall that the Riccati equations are of the form
\begin{equation}\label{eq:re}
w'=a_0(z)w^2(z)+a_1(z)w(z)+a_2(z),
\end{equation}
where $a_0\neq 0, a_2\neq 0,$ $a_i$ are meromorphic functions.
Their well-known properties include:
\begin{itemize}
 \item The Riccati equations  reduce to a linear second order equation.
 \item The solutions of Riccati equations do not posses movable critical points.
\item If one  particular solution is known, the Riccati equation reduces to a linear first order equation.
\item If three particular solutions  $w_1$, $w_2$, $w_3$ are known, then the cross-ratio
$$\big(w(z):w_1(z):w_2(z):w_3(z)\big)$$
is constant along any solution $w$.
There exists a rational function $R$, such that $w=R(w_1, w_2, w_3)$.
\end{itemize}

Mihailo Petrovi\'c, in his PhD thesis in 1894 considered the following  rational ODEs:
\begin{equation}\label{eq:genRiccati}
 w'=\dfrac {P(w, z)}{Q(w, z)},
 \end{equation}
where  $P, Q$ are polynomials in $w, z$. He proved the following theorem:

\begin{theorem}[Petrovi\'c, 1894]\label{th:4}
Such an equation can't have more than three single-valued solutions which present essentially distinct transcendent functions.
\end{theorem}

This result of Petrovi\'c caught an immediate attention of his contemporaries and the Theorem  \ref{th:4}  was quoted, for example, in \cite{Pic1}, \cite{Golubev1}, and \cite{Golubev}.

Let us outline a draft of the proof. As the first step, Petrovi\'c proves that the equations \eqref{eq:genRiccati}
could be reduced to the generalized Riccati equations of the form
\begin{equation}\label{eq:genRiccati1}
w'=\dfrac {P_{n+2}(w, z)}{Q_{n}(w, z)},
\end{equation}
 where $P_{n+2}$, $Q_n$ are polynomials in $w, z$, of degree  $n+2, n$ respectively as polynomials in $w$. This transformation can be done by a change of variables. Then he considered  four cases:
\begin{enumerate}
 \item $Q_n$ has more than two distinct roots: he proved that then all single-valued solutions are rational.
 \item  $Q_n$ has exactly two distinct roots: then all single-valued solutions reduce to at most one transcendental function.
\item $Q_n$ has only one root: then all single valued solutions reduce to at most two essentially distinct transcendental {functions.}
\item $Q_n$ does not contain $w$, thus the equation corresponds to the Riccati equations: then all single valued solutions reduce to at most three essentially distinct transcendental {functions.}
\end{enumerate}

For the illustration, let us show how Petrovi\'c treated the first of the above four cases.

Let $f_1(z), f_2(z), f_3(z)$ be the roots of the polynomial $Q_n$ understood as a polynomial in $w$ with a parameter $z$, i.e. the solutions of $Q_n(w,\, z)=0$ in $\mathbb C(z)[w]$,  and
$w$ be a single-valued solution of the differential equation. Then one can consider $\theta(z)=(w(z):f_1(z):f_2(z):f_3(z))$ which is an algebraic function, since it does not have essential singularities by the Big Picard Theorem and  Lemma about critical point (see Lemma \ref{lemma:lemma} above). Thus, $w$ is a rational combination of algebraic functions and single-valued, thus it is a rational function.

Let us mention two variations  of the Petrovi\'c Theorem  \ref{th:4}.

\begin{theorem}[Golubev, 1911]\label{th:Golubev}
If the above equation under the conditions that $P, Q$ are polynomials in $w$ with finitely many isolated singularities in the coefficients, has three rational solutions then every single-valued solution is rational.
\end{theorem}

A far reaching generalization of the Petrovi\'c Theorem \ref{th:4} was obtained by Malmquist in \cite{Malquist}. Using a very subtle analytic arguments coming from Boutroux he managed to get a very elegant conclusion about the first three items of the above considerations.
\begin{theorem}[Malmquist J. 1913]\label{th:Malmquist}
If the above equation (\ref{eq:genRiccati}) is not a Riccati equation, then all its single-valued solutions are rational functions.
\end{theorem}
A similar result was reproved by Yosida in \cite{Yoshida} using then new Nevannlina theory \cite{Nevanlinna0, Nevanlinna}, see also Theorem \ref{th:Yoshida} below.
As a matter of fact, Malmquist originally proved a much deeper result:
\begin{theorem}[Malmquist J. 1913]\label{th:Malmquist1}
If the above equation (\ref{eq:genRiccati}) with $P, Q$ being polynomial in $w$ with rational coefficients in $z$ has at least one non-algebraic solution which is algebraic over the field of meromorphic functions, then it can be transformed to  a Riccati equation \eqref{eq:re} with rational coefficients, by a transformation of the form:
\begin{equation}\label{eq:transformMalmquist}
 v=\dfrac {P_{n}(w, z)}{Q_{n-1}(w, z)},
\end{equation}
 where $P_{n}$, $Q_{n-1}$ are monic polynomials in $w$, of degree  $n, n-1$ respectively with coefficients rational in $z$.
\end{theorem}
Further results of Malmquist are contained in \cite{Malquist1}, \cite{Malquist2}, \cite{Malquist3}.
 Hille \cite{Hille} prepared a very nice modern survey of the field.
While Golubev and Malmquist quoted Petrovi\'c's result,  Hille \cite{Hille} did not  mention that result. For the most recent developments
of this subject see \cite{Kec} and references therein.

\section{On single-valued transcendental solutions of binomial ODEs of the first order}

In Part 1 of his thesis \cite{Petrovich1}, Petrovi\'c studied also the binomial equations
\begin{equation}
(y')^m=\frac{P(x,X,y)}{Q(x,X,y)},
\label{binom}
\end{equation}
where $m\in\mathbb{N}$, $m\geqslant 2$, $P, Q$ are polynomials in  $x$, $X$ and $y$; the variables  $x$ and $X$ are assumed to be connected through an algebraic relation  $G(x,X)=0$.
Consider the case $m=2$.
Then the equation \eqref{binom} can be rewritten in the form
\begin{equation}
y'=\frac{B(x,X,y)\sqrt{R(x,X,y)}}{C(x,X,y)}.
\label{binom1}
\end{equation}

Petrovi\'c proved the following theorems.

\bigskip

 \begin{theorem}\label{th:5}  If in the equation \eqref{binom1} the number of distinct nonconstant functions  $y_i=\varphi_i(x,X)$ which are the roots either of the polynomial $C$  or the polynomial $R$ is greater than two, then all single-valued in $x$ and $X$ solutions of this equation are rational.
\end{theorem}

\bigskip

\begin{theorem}\label{th:6} If in the equation \eqref{binom1} the polynomial $R$ has one or two nonconstant roots, then this equation does not possess transcendental single-valued solutions.
\end{theorem}
\bigskip

\bigskip
\begin{theorem}\label{th:7}
 In order that the equation \eqref{binom1} has single-valued transcendental solutions, it is necessary that the equation has the form
\begin{equation}
y'=\frac{B(x,X,y)\sqrt{\rho(y)}}{(y-\varphi_1)^{k_1}(y-\varphi_2)^{k_2}},
\end{equation}
 where the polynomial  $B$ has the degree $k_1+k_2$ in $y$, and $\rho(y)$ is a degree four polynomial.
 \end{theorem}
\bigskip

Theorems \ref{th:5}--\ref{th:7} were proven analytically. These theorems could be considered as generalizations of Theorem  \ref{th:4}.
Petrovi\'c did not consider the equations \eqref{binom} with $m>2$ in his thesis. The study of equations  \eqref{binom} in the case  $m>2$ are technically more involved. But he observed that statements and proofs in these cases are still similar to Theorems \ref{th:5}--\ref{th:7} and their proofs for the case $m=2$.

\bigskip
Similar results were obtained 38 years later, by Yosida in 1932.

\begin{theorem}[Yosida 1932, \cite{Yoshida}]\label{th:Yoshida} If algebraic ODE of the form
$${ y'^m}=R(x,y),\;m\in\mathbb{N},$$
where $R$ is a polynomial in $y$, has a transcendental meromorphic solution, then the degree of the polynomial $R$ is not greater than $2m$.
\end{theorem}

\bigskip
Yosida in his paper \cite{Yoshida} quoted the work of Malmquist, but he did not mention the dissertation of Petrovi\'c.

\section{About solutions with fixed singular points of binomial ODEs of the first order}

Petrovi\'c in \cite{Petrovich1} also extracted those binomial equations without movable singular points.

\bigskip

\begin{theorem}\label{th:8}   Among all the equations from the class
\begin{equation}
{ y'^m}=R(x,y),\;\;m\in\mathbb{N},
\label{eq6}
\end{equation}
where $R$ is a function  rational in $y$, only linear ODEs and the equations of the following two types
\begin{equation}\label{eq7}
  { y'^m}=\chi(x)(y-a)^{m-1},
\end{equation}
\begin{equation}\label{eq8}
  { y'^2}=\chi(x)(y-a)(y-b),
\end{equation}
$a,$ $b\in\mathbb{C}$, are such that all solutions have fixed singular points only, i.e. exactly these are the equations for which the set of singular points of the equation coincides with the set of singular points of solutions.
\end{theorem}

\bigskip
   In the first part of the proof, Petrovi\'c proves that $R$  has to be a polynomial if \eqref{eq6} has solutions with fixed singularities.  In the case where the degree of $R$ is not less than $m$, Petrovi\'c reduces the equation \eqref{eq6} to a linear equation  $y'=\sqrt[m]{\chi(x)}(y-\eta(x))$. In the case where the degree of $R$ is less than $m$, the equation \eqref{eq6} reduces to
 \begin{equation}\label{eq9}
 \displaystyle \left(\frac{dy}{dz}\right)^m=S(y),
\end{equation}
 where  {$z$ is a multi-valued function with respect to $x$,}  $S(y)$ is a polynomial in $y$. Petrovi\'c then skillfully applied the Hermite Theorem \ref{th:Hermite} and the results of  Briot and Bouquet \cite{BriotBouquet} to narrow down the class of equations {\eqref{eq9}} and get at the end only those with single-valued solutions and fixed singular points.

\section{On singularities of algebraic ODEs of higher orders}\label{sec:higherorder}

We conclude this work with considerations of singular points of higher order algebraic ODEs. Contrary to the first order case, which, as we mentioned above, was quite completely resolved by Petrovi\'c and his predecessors, the higher order case is still widely open even now, more than 120 years later. There are, however some important subcases which were successfully studied and we are going to list some of them below. As we have already said, see  Section \ref{sec:obstructions}, Petrovi\'c was fully aware of the obstacles preventing his polygonal method to produce
complete results in higher orders and he listed them clearly. Nevertheless, the Petrovi\'c polygonal method can be successfully applied to get some partial answers about higher order equations and to consider certain types of singularities. Petrovi\'c observed that his method could be applied to determine the poles of the solutions in the case of equations {not depending explicitly} on independent variables.

In order to motivate the next question posed by Petrovi\'c, let us go back to the first order case and recall that the Weierstrass equation
\begin{equation}\label{eq:we}
{y'^2}=P_3(y),
\end{equation}
where $P_3$ is a degree three polynomial without multiple zeros, does not depend explicitly on independent variable and has the Weierstrass $\wp$-function as the solution.
Probably motivated by the Hermite Theorem \ref{th:Hermite} as well,  Petrovi\'c applied his polygonal method to study elliptic solutions of higher order algebraic equations, not depending explicitly
on independent variables in \cite{Petrovich2}, (for a recent English translation, see \cite{Pet1a}).

Petrovi\'c singled out the following property.

\bigskip

{\bf Property I.} {\it The polygonal line has at least one edge with a negative integer angular coefficient or it has at least a multiple vertex such that the corresponding characteristic equation has one or several negative integer roots, lying  between the values of the angular coefficients of the two edges that form the multiple vertex.}

\bigskip
\begin{theorem}[Petrovi\'c 1899, \cite{Petrovich2}, \cite{Pet1a}]\label{th:Acta1}
If the equation
\begin{equation}\label{eq18}
 Q(y,y',\dots,y^{(n)})=0,
\end{equation}
where
 $Q$ is a polynomial, has an elliptic solution, then its  polygon  has the property I.
\end{theorem}

Petrovi\'c also considered transformations of solutions.
 Let $R$ be a rational function and   $z=R(y, y',$ $\dots, y^{(q)})$. Let $\Psi(z,z',\dots,z^{(q)})=0$ be  a transformation of the equation { \eqref{eq18}}. He  derived the following result.

\begin{theorem}[Petrovi\'c 1899, \cite{Petrovich2}, \cite{Pet1a}] \label{th:Acta2}
If the polygon, corresponding to the equation $\Psi(z,z',\dots,z^{(q)})=0$ does not possess the property I, then the equation
 \begin{equation}\label{eq17}
 R(y,y',\dots,y^{(q)})={\rm const}
\end{equation}
plays the role of a partial first integral along double-periodic solutions of the equation \eqref{eq18}, i.e. all solutions of the equation \eqref{eq18} of a double-periodic nature also satisfy \eqref{eq17}.
\end{theorem}

 These partial first integrals could serve to reduce the order to get eventually  an equation of the form $Q_1(y, y')=0$, and to treat it further along the lines indicated by
 Briot-Bouqet (see \cite{Ince}, part 2, Ch XIII).

\begin{example}[Petrovi\'c 1899, \cite{Petrovich2}, \cite{Pet1a}] As an example, Petrovi\'c considered the equations of the form
\begin{equation}\label{eq:Petrovic1899}
P_m(y'')=Q_n(y),
\end{equation}
where
$P_m, Q_n$ are given polynomials of degrees  $m, n$ respectively.
The polygon $\mathcal{N}$ is a triangle
with vertices $A(0, 0)$, $B(n, 0)$, $C(m, 2m)$, see Figure \ref{fig:fig5}. In order to  satisfy the Property I,
the triangle $\triangle ABC$ has to be acute. The only edge with the negative angular coefficient is  $BC$, provided $n>m$. The angular coefficient is equal to
$$\frac{2m}{m-n}\in \mathbb Z.$$ With $n>m$, the examples of $(m, n)$ such that  $2m/(m-n)\in \mathbb Z$ include $(m, n)\in\{(1, 2), (1, 3), (2, 4), (2, 6)\}$.

\begin{figure}[h]

    \qquad{\includegraphics[width=12cm]{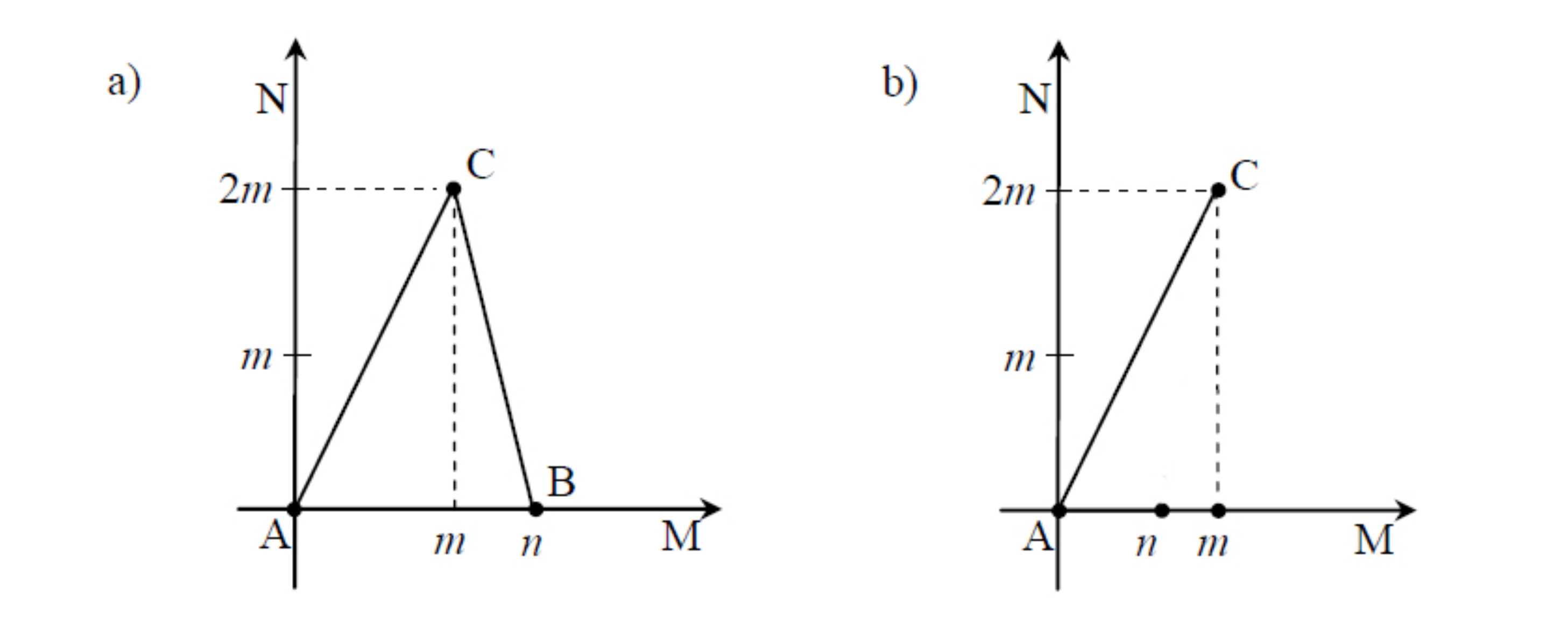}}
    \caption{The Petrovi\'c polygon of the equation \eqref{eq:Petrovic1899}: a) $n>m$; b) $n<m$.}\label{fig:fig5}
\end{figure}
\end{example}

\subsection{Instead of a conclusion}\label{sec:conclusion}

A year after this paper of Petrovi\'c appeared, Painlev\'e published his seminal paper \cite{Pain1} followed by \cite{Pain},
which didn't mark only the turn of the centuries, but opened a new era in the analytic theory of differential equations of higher order.
On the level of ideas, Painlev\'e put forward the Kowalevski position from \cite{Kow} and set the program
to investigate the second order ODEs of the form
$$
y''=Q(x, y, y'),
$$
where $Q$ is a function rational in $y, y'$ and meromorphic in $x$, with the following
property:

$(P)$ -- {\it all solutions are single-valued around all movable singular points}.

This property is called {\it the Painlev\'e property}. Since there were some gaps and errors in Painlev\'e's calculations
in \cite{Pain1}, \cite{Pain}, the program was completed by the students of Painlev\'e, Gambier  in \cite{Gam} and Fuchs (Richard, a son of Lazarus) in \cite {FuchsR}.
As the outcome, they provided 50  equations from which any equation
with the Painlev\'e property  can be obtained by gauge transformations based on M\"obius transformations.
Among these 50 equations, there are six (some of which are multi-parameter families) which can not be solved in terms of existing classical functions (solutions of linear equations or elliptic functions). These equations are
known as {\it the Painlev\'e equations I--VI}. For a general set of parameters the solutions
are new transcendental functions not expressible in terms of formerly known functions,
called {\it the Painlev\'e transcendents}.

The question whether indeed the Painlev\'e equations possess the Painlev\'e property was immediately addressed by Painlev\'e for
the Painlev\'e I equation. His considerations were not quite complete. These questions for all six Painlev\'e equations occupied attention of scientists for more than a century.
For example Golubev in \cite{Golubev2} showed that the solutions of the equations Painlev\'e I-V {possess} the Painlev\'e property. His method was analytic, with a few noncomplete spots, which could be briefly restored. In early 1980's Jimbo and Miwa in \cite{JM} proved the Painlev\'e property for the Painlev\'e VI equations using the connection with the Schlesinger equations. Malgrange obtained a similar result about the same time    in \cite{Malgrange0}.

A full published answer to this question for all six families of Painlev\'e equations appeared in the work of Shimomura \cite{Shim}, see also \cite{HL2004}.

Although a close friend of Painlev\'e and very active in the period of early 1900's, Petrovi\'c did not pay enough attention to the new program of Painlev\'e. Neither him nor his students made any contribution in that direction.

Nevertheless, the ideas of Petrovi\'c polygonal method have been quite recently applied to the Painlev\'e equations in \cite{BrunoGoryuchkina}, \cite{BrunoParusnikova}, and  \cite{BrunoGoryuchkina2}, through
a recently established mathematical discipline,  Power Geometry, as a  modern reincarnation of Fine and Petrovi\'c's method. As we have already mentioned, at that time the authors were unaware of the results of Petrovi\'c and Fine.

For example in \cite{BrunoGoryuchkina2} a method was suggested which allowed to compute elliptic asymptotics of formal solutions of the equations Painlev\'e I--IV. The method is based on the extraction of  approximate equations with the use of the polyhedra associated with systems of ODEs of the first order. The main idea consists in  transformations of the equations Painlev\'e I-IV, which all have a form  $y''=f(x,y,y')$, where $f$ is rational in $y,\;y'$ and meromorphic in  $x$, by use of power transformations of the form
\begin{equation}\label{add0}
y=x^{\alpha}v,\;u=x^{\beta}, \;\beta>0,
\end{equation}
which lead to the equations of the form
 \begin{equation}
 \ddot{v}=h_0(v,\dot{v})+\sum\limits_{i=1}^m u^{-\gamma_i}h_i(v,\dot{v}),\;\;\gamma_i>0.\label{add1}
 \end{equation}
  Here, the choice of numbers  $\alpha$ and  $\beta$ in the transformation   \eqref{add1} is related to the coordinates of the external normals to the faces of the polyhedra
  of the systems of ODEs corresponding to  the equations Painlev\'e I-IV.
  If the approximate equation  $\ddot{v}=h_0(v,\dot{v})$ has a solution $\varphi_0(u)$, which is periodic or double-periodic with a singular point $u=\infty$ (which obviously is also a singular point for the equation \eqref{add1}), then this function is a candidate for the first term of the series
   $$\sum\limits_{j=0}^{\infty}\varphi_j(u)u^{-j},$$
  where $\varphi_j(u)$  are periodic or double - periodic functions with a singular point $u=\infty$, which formally {satisfies} the equation \eqref{add1}.
 It should be mentioned that Boutroux found elliptic asymptotics of the formal solutions of the Painlev\'e I and II equations in \cite{Boutroux}.

Applying the same methods to the Painlev\'e VI equations, Bruno and his  collaborator and student I.~Goryuchkina, found all formal solutions for all values of parameters in neighborhoods of all
singular and nonsingular points of the equations. These formal solutions are not reduced to power series only. Among them, there are generalized power series (power series with complex exponents), Dulac series (integer power series with polynomial in logarithms coefficients), exotic series (integer power series with coefficients, which are meromorphic functions in $x^{i\gamma}$, here $i$ is the imaginary unit and we assume $\gamma\neq 0$), composite series (integer power series with coefficients, which are formal Laurent series with a finite principal  part in powers of $\ln^{-1}{x}$). The formal solutions of the Painlev\'e VI equations of the first three types are uniformly convergent in some open sectors with the vertex in the considered point.  For the Painlev\'e V  equations, all the solutions of {given} types were found, however the formal solutions of  the Painlev\'e V  equations are not exhausted through these types. They require further study and interpretation.

Let us mention that Painlev\'e type program has not yet been completed for the second order ODEs if they are not explicitly resolved
in terms of the second derivative. For higher order ODEs the similar questions are meaningful although very complex. Only some partial
answers are known so far. Chazy was the first to consider the order three case in \cite{Cha}. Very interesting results for the orders four and five
with special forms of $Q$ (for example $Q$ being a  polynomial) were obtained quite recently, see \cite{Cos0}, \cite{Cos1}, \cite{Cos2} and references therein. There are also very recent extensions of the Painlev\'e program to the so-called {\it quasi-Painlev\'e property} which assumes the study of movable algebraic singularities and which has been performed for some  classes of second order algebraic ODEs in \cite{Shim1}, \cite{Shim2}, \cite{FH1},\cite{FH2}, \cite{FH3}, \cite{Kec0}, \cite{Kec}, and references therein.

Along with the Painlev\'e equations, the generalized polygons of Petrovi\'c and Fine have been intensively applied in the studies of solutions of algebraic partial differential equations \cite{Cano3}, \cite{BrunoShadrina}, and also to the studies of solutions of Pfaff systems  \cite{Cano1}, solutions of $q$-difference equations \cite{Cano5}. The polygons have even been used in the proof of Maillet -- Malgrange Theorem  \cite{Cano2}, which provides estimates on the growth of the coefficients of power series, which is an important ingredient
in the selection of summation methods. It is clear that the methods based on polygons of Petrovi\'c and Fine have wide applications and they continue to develop.

We hope we were able to bring to the attention of the specialists in this actively developing field of mathematics and also to a more general audience the gems
almost buried in the past not only to restore the historic justice which these beautiful pioneering results and their outstanding authors deserve but even more -- to propel
these powerful ideas  and put them in the synergy with modern techniques and questions.

{\bf A retrospective of the chronology.} The first ideas to use the Newton -- Puiseux methods in the theory of differential equations probably goes back to
Broit and Bouquet, see \cite{BroitBouquet1}. Henry Fine proposed his polygons in 1889 in \cite{Fine}. Petrovi\'c constructed his polygons in his dissertation
in 1894, \cite{Petrovich1}. He developed further some of his geometric ideas and his polygonal method in his Acta Mathematica paper in 1899 \cite{Petrovich2}. A student of Petrovi\'c, Beri\'c defended his doctoral thesis on the polygonal method at the University of Belgrade in 1912, \cite{Beric}. Starting from 1990's Cano implemented the ideas of Newton-Puiseux polygons to differential equations, see \cite{Cano}. Starting from 1990's Bruno developed
the ideas of Newton-Puiseux polygons in differential equations within the field he named Power Geometry, see \cite{Bruno1}. Cano and his collaborators, as well as Bruno and his collaborators, are still active in the field. The present work is, to the best of our knowledge, the first to consider all these important developments under a uniform perspective.

{\bf Acknowledgements}
The authors are grateful to Renat Gontsov for numerous valuable comments and suggestions. They are also grateful to Joseph Minich and Borislav Gaji\'c for careful reading of the manuscript and their suggestions. The authors would like to thank Professor Susan Friedlander and the Bulletin of the AMS for kindly providing them the opportunity to feature this paper at the cover page of one of the forthcoming issues of BAMS (see \cite{DGBAMS2020}).  The first author would like to thank Professors Bo\v zidar Jovanovi\'c,  Zoran Petri\'c, and   Vera Kova\v cevi\'c-Vuji\v ci\'c, the heads of the departments for Mechanics, for Mathematics, and for Applied Mathematics \& Computer Science for the invitation to deliver a talk about the work of Mihailo Petovi\'c at the first joint meeting of all three colloquia of the Mathematical Institute of the Serbian Academy of Sciences and Arts and to Academicians Stevan Pilipovi\'c and Gradimir Milovanovi\'c, the co-chairs of the public manifestation in Serbia "2018-the Year of Mihailo Petrovi\'c Alas" for the invitation to talk at the round table on May 22nd, 2018, at the Serbian Academy of Sciences and Arts, see \cite{Dra2019}.

\end{document}